\documentclass[11pt]{amsart}
\usepackage{amsxtra}
\usepackage{amssymb}
\theoremstyle{plain}
\newcommand{\cleqn}{\setcounter{equation}{0}}
\newcommand{\clth}{\setcounter{theorem}{0}}
\newcommand {\sectionnew}[1]{\section{#1}\cleqn\clth}
\newcommand{\nn}{\hfill\nonumber}
\newtheorem{theorem}{Theorem}[section]
\newtheorem{lemma}[theorem]{Lemma}
\newtheorem{definition-theorem}[theorem]{Definition-Theorem}
\newtheorem{proposition}[theorem]{Proposition}
\newtheorem{corollary}[theorem]{Corollary}
\newtheorem{definition}[theorem]{Definition}
\newtheorem{example}[theorem]{Example}
\newtheorem{remark}[theorem]{Remark}
\newtheorem{notation}[theorem]{Notation}
\newcommand \bth[1] { \begin{theorem}\label{t#1} }
\newcommand \ble[1] { \begin{lemma}\label{l#1} }

\newcommand \bpr[1] { \begin{proposition}\label{p#1} }
\newcommand \bco[1] { \begin{corollary}\label{c#1} }
\newcommand \bde[1] { \begin{definition}\label{d#1}\rm }
\newcommand \bex[1] { \begin{example}\label{e#1}\rm }
\newcommand \bre[1] { \begin{remark}\label{r#1}\rm }

\newcommand \bnota[1] { \begin{notation}\label{n#1}\rm }
\renewcommand {\eth} { \end{theorem} }
\newcommand {\ele} { \end{lemma} }

\newcommand {\epr} { \end{proposition} }
\newcommand {\eco} { \end{corollary} }
\newcommand {\ede} { \end{definition} }
\newcommand {\eex} { \end{example} }
\newcommand {\ere} { \end{remark} }

\newcommand {\enota} { \end{notation} }
\newcommand \thref[1]{Theorem \ref{t#1}}
\newcommand \leref[1]{Lemma \ref{l#1}}
\newcommand \prref[1]{Proposition \ref{p#1}}
\newcommand \coref[1]{Corollary \ref{c#1}}

\newcommand \lb[1]{\label{#1}}
\def \Cset {{\mathbb C}}
\def \KK {{\mathbb K}}
\def \Zset {{\mathbb Z}}
\def \Nset {{\mathbb N}}
\def \Qset {{\mathbb Q}}

\def \AA  {{\mathcal{A}}}           %mathcal
\def \B  {{\mathcal{B}}}
\def \VV {{\mathcal{V}}}
\def \UU {{\mathcal{U}}}
\def \RR {{\mathcal{R}}}
\def \SS {{\mathcal{S}}}

\def \De {\Delta}   % Greek letters
\def \de {\delta}
\def \al {\alpha}

\def \la {\lambda}

\def \om {\omega}

\def \de {\delta}

\def \sig {\sigma}

\def \sig{\sigma}
\def \mt  {\mapsto}
\def \lha {\leftharpoonup}
\def \rha {\rightharpoonup}

\def \ci  {\circ}

\def \rcor {\rangle}
\def \lcor {\langle}
\def \ol {\overline}
\def \wt {\widetilde}
\def \wh {\widehat}

\def \id { {\mathrm{id}} }

\def \sign { {\mathrm{sign}} }
\def \Ad { {\mathrm{Ad}} }

\def \Lie { {\mathrm{Lie \,}} }

\def \g  {\mathfrak{g}}   % Lie algebra letters

\def \n  {\mathfrak{n}}

\def \b  {\mathfrak{b}}

\def \Mmn {M_{m,n}}
\DeclareMathOperator \wDia  {{ \stackrel{w}{\Diamond} }}
\DeclareMathOperator \Span { {\mathrm{Span}} }

\renewcommand \Im { {\mathrm{Im}} }

\newcommand \Spec { {\mathrm{Spec}} }
\begin{document}
%%%%%%%%%%%%%%%%%%%%%%%%%%%%%%%%%%%%%%%%%%%%%%%%%%%%%%%%%%%%%%%%%%%%%%%%%%%
%%%%%%%%%%%%%%%%%%%%%%    Title    %%%%%%%%%%%%%%%%%%%%%%%%%%%%%%%%%%%%%%%%
\title[Quantized nilpotent Lie algebras]
{Invariant prime ideals in 
quantizations of nilpotent Lie algebras}
\author[Milen Yakimov]{Milen Yakimov}
\address{
Department of Mathematics \\
Louisiana State University \\
Baton Rouge, LA 70803 and
Department of Mathematics \\
University of California \\
Santa Barbara, CA 93106 \\
U.S.A.
}
\email{yakimov@math.lsu.edu}
\date{}
\begin{abstract}
De Concini, Kac and Procesi defined a family of subalgebras $\UU^w_+$
of a quantized universal enveloping algebra $\UU_q(\g)$, associated to 
the elements of the corresponding Weyl group $W$. They are deformations
of the universal enveloping algebras $\UU(\n_+ \cap \Ad_w(\n_-))$
where $\n_\pm$ are the nilradicals of a pair of dual Borel
subalgebras. Based on results of Gorelik and Joseph
and an interpretation of $\UU^w_+$ as quantized algebras 
of functions on Schubert cells, 
we construct explicitly the $H$ invariant prime ideals of each $\UU^w_+$
and show that the corresponding poset is isomorphic to 
$W^{\leq w}$, where $H$ is the group of group-like elements of 
$\UU_q(\g)$. Moreover, for each $H$-prime of $\UU^w_+$ we construct 
a generating set in terms of Demazure modules related to fundamental 
representations. 

Using results of Ramanathan and Kempf we prove similar
theorems for vanishing ideals of closures of torus orbits 
of symplectic leaves of related Poisson structures on Schubert cells 
in flag varieties.
\end{abstract}
\maketitle
%%%%%%%%%%%%%%%%%%%%   Introduction   %%%%%%%%%%%%%%%%%%%%%%%%%%%%%%%%%%%%%%%%
\sectionnew{Introduction}
\lb{intro}
Let $\g$ be a be a split semisimple Lie algebra over a field $\KK$
of characteristic 0. Fix a pair of opposite Borel subalgebras $\b_\pm$ 
with nilradicals $\n_\pm$. Let $q \in \KK$ be transcendental over 
$\Qset$ and $\UU_q(\g)$ be the corresponding quantized universal 
enveloping algebra over $\KK$ with standard generators 
$X^\pm_1, \ldots X^\pm_r$, $K_1^{\pm 1}, \ldots, K_r^{\pm 1}$.

Given an element $w$ of the Weyl group $W$ of $\g$, one defines the nilpotent 
subalgebra $\n_+ \cap \Ad_w(\n_-)$ of $\g$, where $\Ad$ refers to the adjoint 
action. The $q$-analog of $\UU(\n_+ \cap \Ad_w(\n_-))$ is defined 
in a less straightforward way. For each reduced expression 
$w = s_{i_1} \ldots s_{i_n} \in W$ one defines the Lusztig root vectors
\cite{L, CP} $X^\pm_{i_1}$, $T_{i_1}(X^\pm_{i_2})$, $\ldots$, 
$T_{i_1} \ldots T_{i_{n-1}}(X^\pm_{i_n})$, where $T$ denotes the action \cite{L, CP} 
of the braid group of $W$ on $\UU_q(\g)$. De Concini, Kac, and Procesi \cite{DKP} 
proved that the subalgebras of $\UU_q(\g)$ generated by these root vectors 
(in the plus and minus cases) do not 
depend on the choice of a reduced expression of $w$ and studied their representations  
at roots of unity. Denote the De Concini--Kac--Procesi subalgebras of $\UU_q(\g)$ 
corresponding to $w \in W$ by $\UU^w_\pm$.   

In this paper we investigate the set of prime ideals of $\UU^w_-$ invariant 
under the conjugation action of the group $H = \lcor K_1, \ldots, K_r \rcor$ 
of group-like elements of $\UU_q(\g)$. We identify the (finite) poset of 
those ideals ordered under inclusions with a Bruhat interval, obtain an 
explicit 
description of each ideal using Demazure modules, and construct a small 
generating set for each ideal. Special examples of $\UU^w_-$ are
the algebras of quantum matrices $R_q[\Mmn]$. Even in this case the generating 
sets and explicit description of the ideals present new results.

We prove that the algebras $\UU^w_-$ are quotients of the Joseph--Gorelik
quantum Bruhat cell translates \cite{J,G}. The latter are quantizations of the algebras 
of functions on the translated Bruhat cells of the full flag variety
associated to $\g$ with respect to the standard Poisson structure. Along the way 
we obtain a model for $\UU^w_-$ as quantizations of Poisson structures on 
Schubert cells. Our constructions are similar to the De Concini--Procesi \cite{DP} 
interpretation of $\UU^w_-$ as quantum Schubert cells. They constructed an 
isomorphism between a localization of $H \UU^w_-$ and a localization of a 
quotient of the quantized algebra of functions on a Borel subgroup.
We work with a realization of $\UU^w_-$ 
(without $H$ and localization) in terms of Demazure modules.
To be more precise, let $G$ be the split simply 
connected semisimple algebraic group over $\KK$ with Lie algebra $\g$. 
Denote by $B_\pm$ the Borel subgroups of $G$ corresponding to $\b_\pm$. 
It is well known that the coordinate ring of the 
Schubert cell $B_+ w \cdot B_+ \subset G/B_+$ consists of matrix 
coefficients of Demazure $\b_+$-modules, cf. \S \ref{4.6} for details.
We construct a quantum version of this coordinate ring as follows.
Let $P_+$ be the set of dominant weights of $\g$.
Denote by $V(\la)$ the irreducible 
$\UU_q(\g)$-module with highest weight $\la \in P_+$. The Demazure 
module $V_w(\la)$ is the $\UU_q(\b_+)$-module generated by $T_w v_\la$ 
where $v_\la$ is a heighest weight vector of $V(\la)$ and $T_w$ refers 
to the canonical action of the braid group of $W$ on $V(\la)$, 
see \cite{L, CP}. Denote the subalgebra of $\UU_q(\g)$ generated 
by $X^\pm_1, \ldots, X^\pm_r$ by $\UU_\pm$. Identify 
$U^w_+ := U^+ \cap w U^- w^{-1} \cong B_+ w \cdot B_+$
and define the quantized coordinate 
ring $R_q[U^w_+]$ of the Schubert cell 
$B_+ w \cdot B_+$ as the subset of 
$(\UU_+)^*$ consisting of all matrix coefficients  
$c^{w, \la}_\eta(x):= \lcor \eta, x T_w v_\la \rcor$
for $\eta \in V_w(\la)^*$, which is easily seen to be 
a $\KK$-space \S \ref{3.8}. One can make it into a $\KK$-algebra 
by setting
\begin{multline*}
c^{w, \la_1}_{\eta_1} c^{w, \la_2}_{\eta_2}
=
q^{\lcor \la_2, \la_1 - w^{-1} \mu_1 \rcor}
c^{w, \la_1 + \la_2}_{\eta}, \\
\mbox{where} \; \;
\eta = \eta_1 \otimes \eta_2 |_{ \UU_+ (T_w v_{\la_1} \otimes T_w v_{\la_2})}
\in V_w(\la_1 + \la_2)^*
\end{multline*}
for $\eta_1 \in V_w(\la_1)^*$ of weight $\mu_1$ 
and $\eta_2 \in V_w(\la_2)^*$, see \S \ref{3.8} for details.

Recall that to each $w \in W$ one associates a quantum 
$R$-matrix $\RR^w$ which belongs to a certain completion
of $\UU^w_+ \otimes \UU^w_-$, see \S \ref{2.4}.
Our treatment of $\UU^w_-$ rests upon the fact that
he map $\psi_w \colon R_q[U^w_+] \to \UU^w_-$
given by 
\[
\psi_w ( c^{w, \la}_\eta) = (c^{w, \la}_\eta \otimes \id ) (\RR^w)
\]
is an algebra isomorphism and the fact that $R_q[U^w_-]$ is a quotient 
of the Joseph--Gorelik quantum Bruhat cell translates 
(which are quantizations of $w B_- \cdot B_+ \subset G/B_+$). 
Both facts are proved in Sect. 3.

We then use Gorelik's detailed study \cite{G} of the spectra of the 
quantum Bruhat cell translates and Joseph's results \cite{J2} on generating sets 
for ideals of $R_q[G]$ to obtain the following Theorem:

\bth{0.2} Fix $w \in W$. For each $y \in W^{\leq w}$ define
\begin{equation}
\label{Iw}
I_w(y) =
\{ (c^{w,\la}_\eta \otimes \id)(\RR^w) \mid 
\la \in P_+, \eta \in (V_w(\la) \cap \UU^- T_y v_\la)^\perp \}.
\end{equation}
Then:

(a) $I_w(y)$ is an $H$-invariant prime ideal of $\UU^w_-$
and all $H$-invariant prime ideals of $\UU^w_-$ are
of this form.

(b) The correspondence $y \in W^{\leq w} \mapsto I_w(y)$ is an isomorphism
from the poset $W^{\leq w}$ to the poset of $H$ invariant 
prime ideals of $\UU^w_-$ ordered under inclusion; that is
$I_w(y) \subseteq I_w(y')$ for $y, y' \in W^{\leq w}$
if and only if $y \leq y'$.

(c) $I_w(y)$ is generated as a right ideal by
\[
(c^{w, \om_i}_\eta \otimes \id) (\RR^w) \quad 
\mbox{for} \quad
\eta \in (V_w(\om_i) \cap \UU^- T_y v_{\om_i})^\perp, i= 1, \ldots, r,
\]
where $\om_1, \ldots, \om_r$ are the fundamental weights 
of $\g$.
\eth

Assuming only that $q$ is not a root of unity and without restrictions 
on the characteristic of $\KK$ M\'eriaux and Cauchon \cite{MC}
obtained a classification of the $H$-primes of $\UU^w_+$
using Cauchon's deletion procedure \cite{C1}. Such parametrizations were 
previously obtained for quantum matrices \cite{La} by Launois. But even for 
quantum matrices an explicit formula for the ideals $I_w(y)$ of the type
\eqref{Iw} was unknown. 

The poset structure on $H$-primes was known only for quantum matrices 
due to Launois \cite{La} under the same restriction that $q$ is transcendental 
over $\Qset$ and $\KK$ has characteristic 0. 
%In particular, part (b) of \thref{0.2} answers 
%a question of M\'eriaux and Cauchon \cite{MC} under those restrictions.

Generating sets were known only for $3 \times 3$ quantum matrices 
due to Goodearl and Lenagan \cite{GL} for arbitrary $\KK$, $q$ 
not a root of unity. Launois \cite{La0} proved 
that the invariant prime ideals in quantum matrices are generated 
as one sided ideals by quantum minors for the case of $q$ transcendental.
In an independent work Goodearl, 
Launois, and Lenagan \cite{GLL} determine all quantum minors in a given 
invariant prime ideal of $R_q[\Mmn]$ (for an arbitrary  field 
$\KK$, $q$ not a root of unity) and thus construct generating sets
for those ideals (in the case when $\KK$ has characteristic 0 
and $q$ is transcendental over $\Qset$). In Sect. 5 we show that the 
uniform treatment 
of ideal generators for the prime ideals of all algebras 
$\UU^w_-$ from \thref{0.2}, when specialized to quantum matrices
gives explicit generating sets consisting of quantum minors. 
Our generating sets are smaller than those in \cite{GLL}.

Assume that $A$ is an algebra with a rational action of a torus $T$ 
by algebra automorphisms. Goodearl and Letzter \cite{GL} showed 
that under some minor 
conditions $\Spec A$ has a natural stratification into strata indexed 
by $T$-primes of $A$, see also Brown--Goodearl \cite{BG}.
They furthermore proved that each stratum 
can be identified with the spectrum of a Laurent polynomial ring.
Their results apply to iterated skew polynomial rings again under some 
mild hypotheses which are satisfied for $\UU^w_-$ \cite{MC}. This in 
particular 
provides a stratification of $\Spec \, \UU^w_-$ with the property that 
all strata are tori. Such a stratification 
can be also directly obtained from Gorelik's results using \thref{main}.  

In the case when $G$ is a complex simple group, the flag variety 
$G/B_+$ has a natural Poisson structure studied by Brown, Goodearl and 
the author \cite{BGY, GY}. All Schubert cells 
$B_+ w \cdot B_+ \subset G/B_+$ are Poisson submanifolds. This induces
Poisson structures $\pi_w$ on 
$U^w_+ = U^+ \cap w U^- w^{-1} \cong B_+ w \cdot B_+$.
All of them are invariant under the action of the maximal 
torus $T=B_+ \cap B_-$ of $G$. The torus orbits of symplectic leaves 
of $\pi_w$ were described in \cite{BGY,GY}. In Sect. 5 
we show that all results for $\UU^w_-$ have Poisson analogs 
for the vanishing ideals of the torus orbits of leaves of 
$\pi_w$. In particular, one obtains an explicit description 
of these ideals in terms of Demazure modules, as well as 
generating sets for the ideals. Another consequence is that 
we obtain an isomorphism between the poset of $H$-invariant 
prime ideals of $\UU^w_-$ and the underlying poset 
of the stratification of $(U^w_+, \pi_w)$ into 
$T$-orbits of symplectic leaves (with inverted order relation). 
This is the  first step towards realizing the orbit method 
program for the algebras $\UU^w_-$ which would amount
to constructing a homeomorphism between $\Spec \, \UU^w_-$
and the symplectic foliation of $(U^w_+, \pi_w)$. 

To understand the relation between the situation in the Poisson case
and Gorelik's construction one is led to consider 
certain intersections of Schubert cells with respect to
three different flags, see \leref{int} and \prref{leaves}. Recently 
Knutson, Lam and Speyer \cite{KLS} raised the question of finding 
intersections of multiple $(>2)$ Bruhat decompositions with good geometric
properties. They proved that each stratum of Lusztig's 
stratification of Grassmannians can be considered as an 
intersection of Schubert cells with respect to $n$ cyclically 
permuted Borel subgroups and used this to obtain a number of 
results on the geometry of Lusztig's stratification. 
   
We finish the introduction with several notation conventions. 
For a subgroup $B \subset G$ and $g \in G$ we denote by $g \cdot B$
the coset of $g$ in $G/B$ and by $gB$ the product subset of $G$.
For a subvariety $X$ of $Z$ we denote by $Cl_Z(X)$ the Zariski closure 
of $X$ in $Z$. For a subspace $L$ of a vector space $V$ we denote
$L^\perp := \{ \xi \in V^* \mid 
\lcor \xi, v \rcor = 0, \; \forall v \in L \}$.
A submanifold $M$ of a Poisson manifold $(X, \pi)$ will
be called a complete Poisson submanifold if $M$ is 
a union of symplectic leaves of $(X,\pi)$. 
\\ \hfill \\
{\bf Acknowledgements.} I am grateful to Ken Goodearl
for many helpful discussions and to Victor Kac for pointing 
out the results in \cite{DKP} at an early stage of this project
when I was trying to understand representations of quantized universal 
enveloping algebras of nilradicals of parabolic subalgebras.
The author was partially supported by NSF grant DMS-0701107.
I would also like to thank St\'ephane Launois and Laurent Rigal
for comments on the first draft of the paper.
%%%%%%%%%%%%%%%%%%%%%%%%%%%%%%%%%%%%%%%%%%%%%%%%%%%%%%%%%%%%%%%%%%%%%%%%%%%%%%
%%%%%%%%%%%%%%%%%%%%%%%%%%%%%%%%%%%%%%%%%%%%%%%%%%%%%
\sectionnew{The algebras $\UU_q(\g)$, $R_q[G]$, $R^w_0$, and $\UU^w_\pm$}
\lb{qalg}
%%%%%%%%%%%%%%%%
\subsection{}
\label{2.1}
Let $\KK$ be a field of characteristic $0$ and $q \in \KK$ be
transcendental over $\Qset$. Let $\g$ be a split semisimple 
Lie algebra over $\KK$. Denote the rank
of $\g$ by $r$ and its Cartan matrix by
$(c_{ij})$. The quantized universal enveloping algebra 
$\UU_q(\g)$ is the $\KK$-algebra with generators
\[
X^\pm_i, K_i^{\pm 1}, \; i=1, \ldots, r,
\]
subject to the relations
\begin{align*}
&K_i^{-1} K_i = K_i K^{-1}_i = 1, \, K_i K_j = K_j K_i, \, 
K_i X_j^\pm K^{-1}_i = q^{\pm c_{ij}} X_j^\pm,
\\
&X^+_i X^-_j - X^-_j X^+_i = \de_{i,j} \frac{K_i - K^{-1}_i}
{q_i - q^{-1}_i},
\\
&\sum_{k=0}^{1-c_{ij}}
\begin{bmatrix}
1-c_{ij} \\ k
\end{bmatrix}_q
      (X^\pm_i)^k X^\pm_j (X^\pm_i)^{1-c_{ij}-k} = 0, \, i \neq j.
\end{align*}
Here $q_i = q^{d_i}$
for the standard choice of integers $d_i$ for
which the matrix $(d_i c_{ij})$ is symmetric.
Recall that $\UU_q(\g)$ is
a Hopf algebra with comultiplication given by
\begin{align*}
\De(K_i) &= K_i \otimes K_i,
\\
\De(X^+_i) &= X^+_i \otimes K_i + 1 \otimes X^+_i,
\\
\De(X^-_i) &= X^-_i \otimes 1 + K_i \otimes X^-_i,
\end{align*}
antipode and counit given by
\[
S(K_i) = K^{-1}_i, \,
S(X^+_i)= - X^+_i K^{-1}_i, \,
S(X^-_i)= - K_i X^-_i,
\]
and
\[
\epsilon (K_i)=1, \, \epsilon(X^\pm_i)=0.
\]
Here
\[
[n]_q = \frac{q^n - q^{-n}}{q - q^{-1}}, \,
[n]_q ! = [1]_q \ldots [n]_q, \,
\begin{bmatrix}
n \\ m
\end{bmatrix}_{q}
= \frac{[n]_q !}{[m]_q ! [n-m]_q !}.
\]

Denote by $\UU_\pm$ the subalgebras of $\UU_q(\g)$
generated by $\{X^\pm_i\}_{i=1}^r$. Let $H$ be the group 
generated by $\{K_i^{\pm 1}\}_{i=1}^r$.
Set $\UU_q(\b_\pm) = H \UU_\pm$.
%%%%%%%%%%%%%%%%%%%%%%%
\subsection{}
\label{2.2}
Let $P$ and $P_+$ be the sets of all integral and dominant integral 
weights of $\g$. 
The sets of simple roots, simple coroots, and fundamental weights 
of $\g$ will be denoted by $\{\al_i\}_{i=1}^r,$ 
$\{\al_i\spcheck\}_{i=1}^r,$ and $\{\om_i\}_{i=1}^r,$ respectively.
The weight spaces of a $\UU_q(\g)$-module $V$ are defined by
\[
V_\la = \{ v \in V \mid K_i v = q^{ \lcor \la, \al_i\spcheck \rcor} v, \; 
\forall i = 1, \ldots, r \}
\]
for $\la \in P$.
A $\UU_q(\g)$-module is a weight module if it is the 
sum of its weight spaces. The irreducible finite dimensional 
weight $\UU_q(\g)$-modules are parametrized by $P_+$. 
Denote by $V(\la)$ the irreducible module corresponding 
to $\la \in P_+$. For each $\la \in P_+$ fix a highest weight 
vector $v_\la$ of $V(\la)$.

The quantized coordinate ring $R_q[G]$ is the Hopf  
subalgebra of the restricted dual of $\UU_q(\g)$ spanned by 
all matrix entries $c_{\xi, v}^\la,$ 
$\la \in P_+$, $v \in V(\la), \xi \in V(\la)^*$. Thus
$c_{\xi,v}^\la(x) = \lcor \xi, x v \rcor$ for $x \in \UU_q(\g)$.

We have the canonical left and right actions of $\UU_q(\g)$ 
on $R_q[G]$:
\begin{equation}
\label{action}
x \rha c = \sum c_{(2)}(x)c_{(1)}, \; 
c \lha x = \sum c_{(1)}(x)c_{(2)}, \; 
x \in \UU_q(\g), c \in R_q[G].
\end{equation} 
    
The subalgebra of $R_q[G]$ invariant under the left action of $\UU_+$
will be denoted by $R^+$. It is spanned by 
all matrix entries $c_{\xi , v_\la}^\la$ where $\la \in P_+$, 
$\xi \in L(\la)^*$ and $v_\la$ is the fixed highest weight vector 
of $V(\la)$.
%%%%%%%%%%%%%%%%%%%
\subsection{}
\label{2.3}
Denote the Weyl and braid groups of $\g$ by $W$ and $\B_\g$, 
respectively. There is a natural action of $\B_\g$ on the
modules $V(\la)$, see \cite{CP,L} for details.
The standard generators $T_1, \ldots, T_l$
of $\B_\g$ act by
\[
T_i = \sum_{a, b, c \in \Nset} 
(-1)^b q_i^{ac - b} (X^+_i)^{(a)} (X^-_i)^{(b)} (X^+_i)^{(c)} 
\]
where
\[
(X^\pm_i)^{(n)}= \frac{(X^\pm_i)^n}{[n]_{q_i} !} \cdot
\]
Similarly $\B_\g$ acts on $\UU_q(\g)$. Its generators   
act by
\begin{align*}
& T_i(X_i^+) = - X_i^- K_i, \; T_i(X_i^-) = - K_i^{-1} X_i^+, \; 
T_i(K_j) = K_j K_i^{-c_{ij}}, \\
& T_i(X_j^+) = \sum_{r=0}^{- c_{ij}} (- q_i)^{-r} (X_i^+)^{(-c_{ij}-r)}
X_j^+ (X_i^+)^{(r)}, \; j \neq i, \\
& T_i(X_j^-) = \sum_{r=0}^{- c_{ij}} (- q_i)^r (X_i^-)^{(r)}
X_j^- (X_i^-)^{(-c_{ij}-r)}, \; j \neq i.
\end{align*}

The action of the braid group $\B_\g$ has the properties
\begin{equation}
\label{Tcommut}
T_w V(\la)_\mu = V(\la)_{w(\mu)}, \; \; 
T_w ( x . v ) = (T_w x) . (T_w v),
\end{equation}
for all $w \in W$, $x \in \UU_q(\g)$, $v \in V(\la)$, $\mu \in P$.
%%%%%%%%%%%%%%%%%%%%%
\subsection{}
\label{2.4}
Fix $w \in W$. For a reduced expression
\begin{equation}
\label{wdecomp}
w = s_{i_1} \ldots s_{i_k}
\end{equation}
define the roots
\begin{equation}
\label{beta}
\beta_1 = \al_{i_1}, \beta_2 = s_{i_1} \al_{i_2}, 
\ldots, \beta_k = s_{i_1} \ldots s_{i_{k-1}} \al_{i_k}
\end{equation}
and the root vectors
\begin{equation}
X^{\pm}_{\beta_1} = X^{\pm}_{i_1}, 
X^{\pm}_{\beta_2} = T_{s_{i_1}} X^\pm_{i_2}, 
\ldots, X^\pm_{\beta_k} = T_{s_{i_1} \ldots s_{i_{k-1}}} X^{\pm}_{i_k},
\label{rootv}
\end{equation}
see \cite{L} for details.
Following \cite{DKP}, define the 
subalgebra $\UU_\pm^w$ of $\UU_\pm$ generated by 
$X^{\pm}_{\beta_j}$, $j=1, \ldots, k$. 

\bth{DKP} (De Concini, Kac, Procesi) \cite[Proposition 2.2]{DKP}
The definition of the algebra $\UU_\pm^w$
does not depend on the choice of a reduced decomposition of $w$. 
The algebra $\UU_\pm^w$ has the PBW basis
\[
(X^\pm_{\beta_k})^{n_k} \ldots (X^\pm_{\beta_1})^{n_1}, \; \; 
n_1, \ldots, n_k \in \Nset.
\]
\eth

Set
\[
\exp_{q_i} = \sum_{n=0}^\infty q_i^{n(n+1)/2} 
\frac{n^k}{[n]_{q_i}!} \cdot
\]
The universal $R$-matrix associated to $w$ is given by
\begin{equation}
\RR^w = \prod_{j= k, \ldots, 1} \exp_{q_{i_j}}
\left( (1-q_{i_j})^{-2} 
X^+_{\beta_j} \otimes X^-_{\beta_j} \right)
\label{Rw}
\end{equation}
in terms of the reduced decomposition \eqref{wdecomp} and the 
root vectors \eqref{rootv}, cf. \cite{CP} for details.
In \eqref{Rw}
the terms are multiplied in the order $j = k, \ldots, 1$.
The $R$-matrix $\RR^w$ belongs to $\UU^w_+ \wh{\otimes} \UU^w_-$,
the completion of $\UU^w_+ \otimes \UU^w_-$ with respect to the 
descending filtration \cite[\S 4.1.1]{L}. It does not depend on the 
choice of reduced decomposition of $w$. For all $\la, \mu \in P_+$
\begin{equation}
\lb{RwT}
T_{w, V(\la) \otimes V(\mu)} = (\RR^w)^{-1} 
\left( T_{w, V(\la)} \otimes T_{w, V(\mu)} \right).
\end{equation}  
%%%%%%%%%%%%%%%%%%%%%%%%%%%%%%%%%%%%%%%%%%%%%%%%%%%%%%%%%%%%%%%%%
\subsection{}
\label{2.5}
Fix again $w \in W$. For each $\la \in P_+$ fix 
$\xi_{w, \la} \in (V(\la)^*)_{- w\la}$ normalized by
$\lcor \xi_{w, \la}, T_w v_\la \rcor =1$.

Following \cite{J}, for $\la \in P_+$ define
\[
c^\la_w = c^\la_{\xi_{w,\la}, v_\la},
\] 
cf. \S \ref{2.2} for the definition of $v_\la \in V(\la)_\la$.

It is clear from \eqref{RwT} that 
$c^\la_w c^\mu_w = c^{\la+ \mu}_w= c^\mu_w c^\la_w$
for all $\la, \mu \in P_+$. Denote \cite[\S 9.1.10]{J} the multiplicative 
commutative subset of $R_+$:
\[
c_w = \{ c_w^\la \mid \la \in P_+ \}.
\]

\ble{Ore} (Joseph) \cite[Lemma 9.1.10]{J} 
The set $c_w$ is Ore in $R^+$.
\ele

Denote the localization
\[
R^w = R^+[c_w^{-1}],
\]
cf. \cite{J, G} for details. Recall that the 
left and right $\UU_q(\g)$-actions \eqref{action}
on $R_q[G]$ induce left and right actions of 
$\UU_q(\g)$ on $R^w$, \cite[\S 4.3.12]{J}.
Let $R_0^w$ be the $H$ invariant 
subalgebra of $R^w$ with respect to the left action of $H$. 
For $\la \in P_+$ set $c_w^{-\la} = (c_w^\la)^{-1}$
in $R^w$. Note that
\begin{equation}
\label{Rw0}
R^w_0 = \{ c_w^{-\la} c^\la_{\xi, v_\la} \mid 
\la \in P_+, \xi \in V(\la)^* \}
\end{equation}
since for all $\la, \mu \in P_+$ and 
$\xi \in V(\la)^*$, there exists $\xi' \in V(\la+\mu)^*$ 
such that $c^{-\la}_w c^\la_{\xi, v_\la} =
c^{-\la-\mu}_w c^{\la+\mu}_{\xi', v_{\la+\mu} }$.
%%%%%%%%%%%%%%%%%%%%%%%%%%%%%%%%%%%%%%%%%%%%%%%%%%%%%
\sectionnew{The $H$-spectrum of $\UU^w_-$} 
\label{hsp}
%%%%%%%%%
\subsection{} 
\label{3.1}
We start by recalling several results of Gorelik \cite{G}.
For $y \in W$ define the ideals
\[
Q(y)^\pm = \Span \{ c^\la_{\xi, v_\la} \mid \xi \in V(\la)^*, \, \xi 
\perp \UU_\pm T_y v_\la \}
\]
of $R^+$ and 
\begin{equation}
\label{id2}
Q(y)^\pm_w = \{c^{-\la}_w c^\la_{\xi, v_\la} \mid \xi \in V(\la)^*, \,
\xi \perp \UU_\pm T_y v_\la \}
\end{equation}
of $R^w_0$. In the second case one does not need to take span 
because of \eqref{Rw0}. The ideals \eqref{id2} are nontrivial if and only if 
$y \geq w$ in the plus case and $y \leq w$ in the minus 
case.

Following \cite{G} denote
\[
W \wDia W = \{ (y_-, y_+) \in W \times W \mid y_- \leq w \leq y_+ \}.
\]
Consider the induced poset structure on $W \wDia W$ from the standard
Bruhat order on the first copy of $W$ and the inverse Bruhat order 
on the second copy of $W$, i.e. for 
$(y_-, y_+), (y'_-, y'_+) \in W \times W$
\[
(y_-, y_+) \leq (y'_-, y'_+) \; \; \mbox{if and only if} \; \; 
y_- \leq y'_- \; \mbox{and} \; y_+ \geq y'_+. 
\]

\bth{Gorelik} (Gorelik) \cite[Lemma 6.6, Proposition 6.11, Corollary 7.1.2]{G} 

(a) For each $(y_-, y_+) \in W \wDia W$ there exists a unique 
$H$ invariant prime ideal $Q(y_-, y_+)_w$ of $R^w_0$ 
which is minimal among all $H$ invariant prime ideals of 
$R^w_0$ containing $Q(y_-)^-_w + Q(y_+)^+_w$. 

(b) The ideals in (a) are distinct
and exhaust all $H$ invariant prime ideals 
of $R^w_0$.

(c) The map $Q(y_-,y_+)_w \mapsto (y_-, y_+)$ is an isomorphism
between the posets of $H$ invariant prime ideals 
of $R_0^w$ ordered under inclusion 
and $W \wDia W$; that is $Q(y_-, y_+)_w \subseteq Q(y'_-, y'_+)_w$ for 
$(y_-, y_+), (y'_-, y'_+) \in W \wDia W$ if and only if
$y_- \leq y'_- \leq w \leq y'_+ \leq y_+$.

(d) For $y_- \in W^{\leq w}$ 
\[
Q(y_-, w)_w = Q(y_-)^-_w + Q(w)^+_w.
\]
\eth
%%%%%%%%%%%%%%%%%%%%
\subsection{}
\label{3.2}
We will construct a surjective homomorphism from 
$R_0^w$ to $\UU^w_-$ similarly to the De Concini--Procesi 
construction of quantum Schubert cells \cite[Sect. 3]{DP}. 
We start with the following simple Lemma. 
Its proof is included for completeness. 
\ble{HA} Assume that $H$ is a Hopf algebra and $A$ is an 
$H$-module algebra with a right $H$ action.
Let $\epsilon \colon A \to \KK$ be an algebra homomorphism, 
where $\KK$ is the ground field.
Then the map $\phi \colon A \to H^*$ given by
\[
\phi(a) (h) = \epsilon (a . h)
\]
is an algebra homomorphism. If, in addition the action of $H$ 
is locally finite, then the image of $\phi$ is contained in
the restricted dual $H^\circ$ of $H$.
\ele 
\begin{proof}
Let $a_1, a_2 \in A$. Using Sweedler's notation 
\begin{multline}
\phi(a_1 a_2) (h) = \epsilon ( (a_1 a_2) . h) =
\epsilon( \sum (a_1 . h_{(1)}) (a_2 . h_{(2)}) ) \\ 
= \sum \epsilon (a_1 . h_{(1)}) \epsilon (a_2 . h_{(2)})
= \lcor \phi(a_1) \otimes \phi (a_2), \De(h) \rcor = 
(\phi(a_1) \phi(a_2))(h). 
\end{multline}

If the action of $H$ is locally finite, then 
for all $a \in A$ the annihilator in $H$ of the 
finite dimensional module $a.H$ is an ideal $J$ 
of finite codimension. Since $J \subseteq \ker
\phi(a)$, $\phi(a) \in H^\circ$.
\end{proof}
%%%%%%%%%%%%%%%%%%%%
\subsection{}
\label{3.3}
\ble{ep} The map $\epsilon_w \colon R^+ \to \KK$
given by 
\begin{equation}
\label{epw}
\epsilon_w ( c^\la_{\xi, v_\la}) = \xi(T_w v_\la), \; \; 
\la \in P_+, \xi \in V(\la)^* 
\end{equation}
is an algebra homomorphism. Moreover $\epsilon_w(c_w^\la) = 1$ 
and $\epsilon_w$ induces a homomorphism from 
$R^w$ to $\KK$ (which will be denoted by the same 
letter).
\ele

The fact that \eqref{epw} defines a homomorphism
is straightforward from \eqref{RwT}. The equality
$\epsilon_w(c_w^\la) = 1$ follows from the definition
of $c_w^\la$.

Consider the right action of $\UU_q(\b_+)$ on $R_0^w$, 
cf. \S \ref{2.5} and \cite[\S 4.3.12]{J}.
\leref{HA} and \leref{ep} now imply that the map 
\begin{equation}
\label{phiw}
\phi_w \colon R^w_0 \to (\UU_q(\b_+))^*, \quad
\lcor \phi_w(c), x \rcor = \epsilon_w (c \lha x), \; 
c \in R^w_0, x \in \UU_q(\b_+) 
\end{equation}
is an algebra homomorphism.

\subsection{}
\label{3.4}
Let $w_\circ$ be the longest element of the Weyl group $W$ of $\g$.
Fix a reduced expression for $w \in W$ as in \eqref{wdecomp}. 
There exists a reduced expression of $w_\circ$ starting with 
the expression \eqref{wdecomp}:
\begin{equation}
\label{w0expr}
w_\circ= s_{i_1} \ldots s_{i_k} s_{i_{k+1}} \ldots s_{i_N}.
\end{equation}
Define the roots
\begin{equation}
\label{beta2}
\beta_1 = \al_{i_1}, \beta_2 = s_{i_1} \al_{i_2}, 
\ldots, \beta_N = s_{i_1} \ldots s_{i_{N-1}} \al_{i_N}
\end{equation}
of $\g$ and the root vectors
\begin{equation}
X^{\pm}_{\beta_1} = X^{\pm}_{\al_{i_1}}, 
X^{\pm}_{\beta_2} = T_{s_{i_1}} X^\pm_{\al_{i_2}}, 
\ldots, X^\pm_{\beta_N} = T_{s_{i_1} \ldots s_{i_{N-1}}} X^{\pm}_{\al_{i_N}}.
\label{rootv2}
\end{equation}
of $\UU_q(\g)$. For $j \leq k$, $\beta_j$ and $X^\pm_{\beta_j}$ are exactly 
the roots and the root vectors defined by \eqref{beta} and \eqref{rootv}.

\ble{eval} For all $K \in H$, $j \in \Nset$, $k < j \leq N$, 
$Y \in \UU_q(\b_+)$, $c \in R_0^w$
\[
\lcor \phi_w(c), Y K \rcor = \lcor \phi_w(c),  Y \rcor
\]  
and
\[
\lcor \phi_w(c), Y X^+_{\beta_j} \rcor = 0.
\]  
\ele
\begin{proof} The first equality is simply the definition 
of $R_0^w$ as an invariant subalgebra of $R^w$.

Applying \cite[\S 4.3.12]{J} one sees that the second equality follows
from
\[
\lcor \xi, YX_{\beta_j}^+ T_w v_\la \rcor = 0, \; 
\forall j \geq k+1, Y \in \UU_+, \xi \in V(\la)^*.
\]  
This in turn is proved by observing that 
$T^{-1}_w(X_{\beta_j}^+)= T_{s_{i_{k+1}} \ldots s_{i_{j-1}}} X^+_{i_j} 
\in \UU_+$ annihilates $v_\la$.  
\end{proof}

\bco{eval2} For all $n_1, \ldots, n_N, n \in \Nset$ 
\begin{multline}
\lcor \phi_w( c_w^{-\la} c^\la_{\xi, v_\la}), 
(X_{\beta_1}^+)^{n_1} \ldots (X_{\beta_N}^+)^{n_N} K^n) \rcor
\\
=\lcor \xi,
(X_{\beta_1}^+)^{n_1} \ldots (X_{\beta_N}^+)^{n_N} T_w v_\la \rcor
= \delta_{n_{k+1}, \ldots, n_N, 0} 
\lcor \xi,
(X_{\beta_1}^+)^{n_1} \ldots (X_{\beta_k}^+)^{n_k} T_w v_\la \rcor.
\end{multline}
\eco

\subsection{}
\label{3.5}
The standard bilinear form $\UU_q(\b_+) \times \UU_q(\b_-) \to \KK$ 
can be used to embed $\UU_q(\b_-)$ in $(\UU_q(\b_+))^*$ (as algebras). 
We identify $\UU^w_-$ with its image in $(\UU_q(\b_+))^*$.    

\bpr{main1} The image of $\phi_w \colon R_0^w \to (\UU_q(\b_+))^*$
is $\UU^w_-$. Its kernel is $Q(w)^+_w$.
\epr
\begin{proof} The inclusion $\Im \phi_w \subseteq \UU^w_-$ follows 
from the fact that 
\[
\{(X_{\beta_1}^+)^{n_1} \ldots (X_{\beta_N}^+)^{n_N} \mid n_1, 
\ldots, n_n \in \Nset\}
\]
is a PBW basis of $\UU_+$ and \coref{eval2}.
Assume that $\phi_w$ is not surjective. Then there exists 
$X \in \UU_w^+$, $X \neq 0$ such that $\xi(X T_w v_\la) = 0$ for 
all $\la \in P_+$ and $\xi \in V(\la)^*$. Therefore 
$X_1 = T_w^{-1}(X) \in \UU^-$ satisfies $c(X_1) = 0$ for all 
$c \in R^+$, i.e. $X_1 v_\la =0$ for all $\la \in P_+$. 
It is well known that this implies $X_1=0$, see e.g. \cite[\S 4.3.5-4.3.6]{J},
which contradicts with $X \neq 0$.

Let $\la \in P_+$, $\xi \in V(\la)^*$. Using \coref{eval2} 
we see that
$c^{-\la}_w c^\la_{\xi, v_\la} \in \ker \phi_w$ if and only if 
\[
\lcor \xi, \UU_+ T_w v_\la \rcor = 0,
\]
i.e. $c^{-\la}_w c^\la_{\xi, v_\la} \in Q(w)_w^+$. 
Thus $\ker \phi_w = Q(w)_w^+$.
\end{proof}
%%%%%%%%%%%%%%%%%%%%
\subsection{}
\label{3.6}
\prref{main1} and \S 3.2 imply the following explicit form of the 
map $\phi_w$, cf. Theorem 3.2 of De Concini and Procesi \cite{DP}:

\bth{main} The $\KK$-linear map 
\begin{equation}
\label{phiw2}
\phi_w \colon R_0^w \to \UU^w_-, \; \; 
\phi_w(c_w^{-\la} c^\la_{\xi, v_\la}) = 
(c^\la_{\xi, T_w v_\la} \otimes \id) (\RR^w), \; 
\la \in P_+, \xi \in V(\la)^* 
\end{equation}
is a surjective homomorhism of algebras. Its kernel is $Q(w)^+_w$.
\eth

The explicit form of $\phi_w$ originally defined by \eqref{phiw} 
follows from the fact that the $R$-matrix $\RR^w$ is by definition 
equal to a sum of the form $\sum_i Y_i \otimes Z_i$ where $\{Y_i\}$ and $\{Z_i\}$
are dual bases of $\UU^w_+$ and $\UU^w_-$ with respect
the standard bilinear form $\UU_q(\b_+) \times \UU_q(\b_-) \to \KK$.
%%%%%%%%%%%%%%%%%%%%
\subsection{}
\label{3.7}
\coref{eval2} implies that $\phi_w \colon R_0^w \to \UU^w_-$ 
is $H$-equivariant with respect to the right action of $H$ on 
$R^w_0$ \eqref{action} and the conjugation action of 
$H$ on $\UU^w_-$:
\begin{equation}
\label{Hconj}
K. x = K^{-1} x K, \quad K \in H, x \in \UU^w_-.
\end{equation}

The following Theorem describes the poset of $H$-primes 
of $\UU^w_-$.

\bth{main2} Fix $w \in W$. For each $y \in W^{\leq w}$ define
\begin{multline}
\label{P_w}
I_w(y) = \phi_w ( Q(y)^-_w) = 
\{ (c^\la_{\xi, T_w v_\la} \otimes \id)(\RR^w) \mid \\
\la \in P_+, \xi \in V(\la)^*, \xi \perp \UU^- T_y v_\la \}.   
\end{multline}
Then:

(a) $I_w(y)$ is an $H$ invariant prime ideal of $\UU^w_-$.

(b) All $H$ invariant prime ideals of $\UU^w_-$ are 
of this form.

(c) The correspondence $y \in W^{\leq w} \mapsto I_w(y)$ is an isomorphism 
from the poset $W^{\leq w}$ to the poset of $H$ 
invariant prime ideals of $\UU^w_-$ ordered under 
inclusion; that is $I_w(y) \subseteq I_w(y')$ for 
$y, y' \in W^{\leq w}$ if and only if $y \leq y'$.
\eth

\begin{proof} The map $\phi_w$ establishes a bijection 
between prime ideals of $\UU^w_-$ and prime ideals of $R^w_0$ containing
$\ker \phi_w = Q(w)^+_w$ (in order preserving way).
The map $\phi_w$ is also equivariant with respect to the right action 
of $H$ on $R^w_0$ and the conjugation action of $H$ on $\UU^w_-$, cf. 
\ref{3.5}. Thus it provides an isomorphism between the posets of $H$
invariant prime ideals 
of $R^w_0$ containing $\ker \phi_w = Q(w)^+_w$ and the $H$-primes 
of $\UU^w_-$. Now \thref{main2} follows from Gorelik's \thref{Gorelik} 
because the only $H$ invariant prime ideals of $R^w_0$ that contain $Q(w)^+_w$
are $Q(y_-,w)_w$, $y_- \in W^{\geq w}$ (see \thref{Gorelik} (a), (c)).
\end{proof}

\bre{1} Gorelik proved \cite{G} that all ideals $Q(y_-, y_+)_w$
are completely prime and as a consequence one gets that 
all ideals $I_w(y)$ are completely prime. This is true 
in a greater generality for $H$-primes of certain iterated 
skew polynomial rings \cite[Proposition 4.2]{GLe} 
due to Goodearl and Letzter.
\ere

\bre{2} There is a natural action of the algebraic torus 
$\KK^r$ on $\UU_q(\g)$ by algebra automorphisms constructed 
by setting
\[
(a_1, \ldots, a_r) \cdot X_i^\pm = a_i^{\pm 1} X_i^+, 
(a_1, \ldots, a_r) \cdot K_i = K_i,
\quad i= 1, \ldots r.
\]
The subalgebras $\UU^w_\pm$ are invariant under it. A subset of 
$\UU_q(\g)$ is invariant under the action of $\KK^r$ if and only if
it is invariant under the conjugation action of $H$. From the point 
of view of Goodearl--Letzter theory of $H$-primes it is more natural 
to use the action of $\KK^r$ since this group is algebraic. Because 
the invariance properties under $\KK^r$ and $H$ are the same and the 
latter action is more natural within the Hopf algebra setting, we use 
the $H$ action.
\ere
%%%%%%%%%%%%%%%%%%%%%%%%%%%%%%%%%%%
\subsection{}
\label{3.8}
An equivalent way to define the algebras $\UU^w_-$ and to work with them 
is by using Demazure modules. This gives an interpretation
of $\UU^w_-$ as quantized algebras of functions on Schubert cells
which is similar to the De Concini--Procesi isomorphism 
\cite[Theorem 3.2]{DP}. A notion of quantum Schubert cells 
in the case of Grassmannians was also defined in 
\cite{LLR} using the algebra of quantum matrices. 

Recall that the $\UU_q(\b_+)$-modules 
$V_w(\la) = \UU_+ T_w v_\la = \UU^w_+ T_w v_\la$ are called 
Demazure modules, cf. \cite[\S 4.4 and 6.3]{J} for details.
For $\eta \in V_w(\la)^*$ define 
\[
c^{w,\la}_{\eta} \in (\UU_+)^*, \; \; c^{w,\la}_{\eta}(X) = \lcor \eta, 
X T_w v_\la \rcor, \; X \in \UU_+.
\]
Set $U^w_+ = U_+ \cap w U_- w^{-1}$. Denote 
by $R_q[U^w_+]$ the subset of $(\UU_+)^*$ consisting of 
\[
c^{w, \la}_\eta, \; \; \la \in P_+, \eta \in V_w(\la)^*.
\]
Consider the linear map
\begin{equation}
\label{varphi}
\varphi_w \colon R^w_0 \to (\UU_+)^*, \; \; 
\varphi_w ( c^{-\la}_w c^\la_\xi ) = c^{w, \la}_{\xi |_{V_w(\la)} }, 
\; 
\la \in P_+, \xi \in V(\la)^*.
\end{equation}
(Because of \coref{eval2} this
is nothing but the map $\phi_w \colon R^w_0 \to (\UU_q(\b_+))^*$ 
composed with the the linear projection $(\UU_q(\b_+))^* \to (\UU_+)^*$.)
The image of $\varphi_w$ is $R_q[U^w_+]$ and 
its kernel is $Q(w)^+_w$ because of 
\coref{eval2}. In particular $R_q[U^w_+]$ is a 
subspace of $(\UU_+)^*$. Since $Q(w)_w^+$ 
is an ideal of $R^w_0$ one can push forward 
the algebra structure of $R^w_0$ to an algebra 
structure on $R_q[U^w_+]$. 
From now on $R_q[U^w_+]$ will denote the 
subspace of $(\UU_+)^*$ equipped with this algebra 
structure. Recall \cite[\S 6.5]{G} that for all 
$\la_1, \la_2 \in P_+$, $\xi_1 \in (V(\la_1)^*)_{\mu_1}$,
$\xi_2 \in V(\la_2)^*$
\begin{equation}
\label{ccomm}
c^{-\la_1}_w c_{\xi_1, v_{\la_1} }^{\la_1} 
c^{-\la_2}_w c_{\xi_2, v_{\la_2} }^{\la_2}
= q^{\lcor \la_2, \la_1 - w^{-1} \mu_1 \rcor}
c^{-\la_1 -\la_2 }_w c_{\xi_1, v_{\la_1} }^{\la_1} 
c_{\xi_2, v_{\la_2} }^{\la_2}, 
\end{equation}
see \cite[\S 9.1]{J} for more details on commutation relations in $R_q[G]$.
Let $\eta_1 \in (V_w(\la_1)^*)_{\mu_1}$ and $\eta_2 \in V_w(\la_2)^*$.
The induced algebra structure on $R_q[U^w_+]$ is given by:
\begin{multline}
\label{Uwmult}
c^{w, \la_1}_{\eta_1} c^{w, \la_2}_{\eta_2} 
= 
q^{\lcor \la_2, \la_1 - w^{-1} \mu_1 \rcor}
c^{w, \la_1 + \la_2}_{\eta}, \\
\mbox{where} \; \; 
\eta = \eta_1 \otimes \eta_2 |_{ \UU_+ (T_w v_{\la_1} \otimes T_w v_{\la_2})}
\in V_w(\la_1 + \la_2)^*.
\end{multline}
Here we use that $\UU_+ (T_w v_{\la_1} \otimes T_w v_{\la_2}) 
\subset \UU_+ T_w v_{\la_1} \otimes \UU_+ T_w v_{\la_2}$.
Note that it is not a priori obvious from \eqref{Uwmult} that this is a well defined 
multiplication and that it is associative. The unit in 
$R_q[U^w_+]$ is equal to $c^{w, \la}_{\eta_{w,\la}}$ 
where $\eta_{w, \la} = \xi_{w, \la}|_{V_w(\la)}$, 
$\la \in P_+$, recall \S 2.5. (The elements 
$c^{w, \la}_{\eta_{w,\la}} \in (\UU_+)^*$ are all equal to each other.)

Define the linear map 
\begin{equation}
\label{psi}
\psi_w \colon R_q[U^w_+] \to \UU^w_-, \quad
\psi_w ( c^{w, \la}_\eta ) =
(c^{w, \la}_\eta \otimes \id) (\RR^w), \; 
\la \in P_+, \eta \in V_w(\la)^*.
\end{equation} 
\coref{eval2} and \thref{main}
imply that $\psi_w \colon R_q[U^w_+] \to \UU^w_-$ is an 
algebra isomorphism.

Recall that $\phi_w \colon R^w_0 \to \UU^w_-$
is $H$-equivariant with respect to the right action
of $H$ on $R^w_0$ \eqref{action} and the conjugation action 
of $H$ on $\UU^w_-$ \eqref{Hconj}. 

From the definition of $\varphi_w$ one obtains 
that $\varphi_w \colon R^w_0 \to R_q[U^w_+]$
is $H$-equivariant with respect to the right action 
of $H$ on $R^w_0$ and the restriction to $R_q[U^w_+]$ of
the following action of $H$ on $\UU_+$ 
\begin{equation}
\label{HactR}
K.c = K^{-1} \rha c \lha K, 
\lcor K^{-1} \rha c \lha K, X \rcor = c ( K X K^{-1} )
\end{equation} 
$K \in H, c \in (\UU_+)^*, X \in \UU_+$.
Finally, the isomorphism $\psi_w \colon R_q[U^w_+] \to \UU^w_-$
is $H$-equivariant with respect to the action \eqref{HactR}
and the conjugation action of $H$ on $\UU^w_-$ \eqref{Hconj}.

\bth{main3} (1) The homomorphism 
$\phi_w \colon R^w_0 \to \UU^w_-$ factors through
the surjective homomorphism \eqref{varphi}
$\varphi_w \colon R^w_0 \to R_q[U^w_+]$ and 
the isomorphism \eqref{psi} 
$\psi_w \colon R_q[U^w_+] \to \UU^w_-$.
Both maps are $H$-equivariant.

(2) Under the isomorphism 
$\psi_w \colon R_q[U^w_+] \to \UU^w_-$ 
the ideals $I_w(y)$, $y \in W^{\leq w}$
of $\UU^w_-$ correspond to the $H$
invariant prime ideals 
\[
J_w(y) = \{c^{w,\la}_\eta \mid \la \in P_+, 
\eta \in (V_w(\la) \cap \UU_- T_y v_\la)^\perp
\}
\]
of $R_q[U^w_+]$.

(3) The ideals $I_w(y)$, $y \in W^{\leq w}$
of $\UU^w_-$ are also given by
\[
I_w(y) = \{ (c^{w,\la}_\eta \otimes \id)(\RR^w) \mid 
\la \in P_+, \eta \in (V_w(\la) \cap \UU_- T_y v_\la)^\perp 
\}. 
\]
\eth
Part (1) has already been established. Recall that 
\[
I_w(y) = \phi_w (Q(y)_w^-+ Q(w)_w^+). 
\]
Part (2) is a direct computation of $J_w(y):= \varphi_w (Q(y)_w^-)$. 
Part (3) follows from the fact that $I_w(y) = \psi_w(J_w(y))$.   
%%%%%%%%%%%%%%%%%%%%%%%%%%%%%%%%%%
\subsection{}
\label{3.9} 
In this subsection, based on results of Joseph \cite{J2}, we
construct generating sets of the $H$-primes $I_w(y)$ of $\UU^w_-$.
Denote the subalgebra
\[
R^- = \{ c^\la_{\xi, T_{w_\circ} v_\la} \mid \la \in P_+, \xi \in V(\la)^* \}
\]
of $R_q[G]$. Recall \S \ref{2.2} that  $\om_1, \ldots, \om_l$ denote 
the fundamental weights of $\g$.

\bth{Joseph} (Joseph) \cite[Th\'eor\`eme 3]{J2} For all $w \in W$ 
\begin{equation}
\label{Jeq1}
\Span \{ c^\la_{\xi, v_\la} \mid 
\la \in P_+, \xi \in (\UU_+ T_w v_\la)^\perp \}
\\ = \sum_{i=1}^r \{ c^\la_{\xi, v_{\om_i}} \mid 
\xi \in (\UU_+ T_w v_{\om_i})^\perp \} R^+
\end{equation}
and
\begin{multline}
\label{Jeq2}
\Span \{ c^\la_{\xi, T_{w_\circ} v_\la} \mid 
\la \in P_+, \xi \in (\UU_- T_w v_\la)^\perp \}
\\ = \sum_{i=1}^r \{ c^\la_{\xi, T_{w_\circ} v_{\om_i}} \mid 
\xi \in (\UU_- T_w v_{\om_i})^\perp \} R^-.
\end{multline}
\eth

The left hand sides of \eqref{Jeq1}--\eqref{Jeq2} 
are $H$ invariant prime ideals of $R^\pm$ and the right hand 
sides give efficient generating sets of them as right ideals.

Using \eqref{RwT}, one sees that the map 
\[
c^\la_{\xi, T_{w_\circ} v_\la} \in R^- \mapsto
c^\la_{\xi, v_\la} \in R^+
\]
is an isomorphism of algebras.
Thus \eqref{Jeq2} implies 
\[
\Span \{ c^\la_{\xi, v_\la} \mid 
\la \in P_+, \xi \in (\UU_- T_w v_\la)^\perp \}
= \sum_{i=1}^r \{ c^\la_{\xi, v_{\om_i}} \mid 
\xi \in (\UU_- T_w v_{\om_i})^\perp \} R^+.
\]
Therefore for all $y \in W^{\leq w}$:
\begin{multline*}
Q(y)_w^- = \{ c^\la_{\xi, v_\la} c^{-\la}_w \mid
\la \in P_+, \xi \in (\UU_- T_y v_\la)^\perp \} \\
\subset \sum_{i=1}^r \{ c^{\om_i}_{\xi, v_{\om_i}} c^{-\om_i}_w \mid 
\xi \in (\UU_- T_y v_{\om_i})^\perp \} R^+[c^{-1}_w] 
\end{multline*}
Using the left action of $H$ on 
$R^+[c^{-1}_w]$ \eqref{action}
and the fact that $R^+[c^{-1}_w]$ is a semisimple $H$-module
we obtain 
\[
Q(y)_w^- =
\sum_{i=1}^r \{ c^{\om_i}_{\xi, v_{\om_i}} c^{-\om_i}_w \mid 
\xi \in (\UU_- T_y v_{\om_i})^\perp \} R^w_0.
\]

Using \eqref{ccomm} we see that the ideals $I_w(y)$ of 
$\UU^w_-$ (and $J_w(y)$ of $R_q[U^w_+]$)
are generated as right ideals by the subsets in 
\thref{main2} and \thref{main3} corresponding 
to the fundamental weights $\om_1, \ldots, \om_l$.

\bth{main4} For all $w \in W$ and $y \in W^{\leq w}$:
\[
I_w(y) = \sum_{i=1}^r 
\{ (c^{\om_i}_{\xi, T_w v_{\om_i}} \otimes \id) (\RR^w) \mid 
\xi \in (\UU^- T_y v_{\om_i})^\perp \} \UU^w_-
\]
and 
\[
J_w(y) = \sum_{i=1}^r \{ c^{w,\om_i}_\eta \mid 
\eta \in (V_w(\om_i) \cap \UU_- T_y v_{\om_i})^\perp\} R_q[U^w_+].
\]
\eth
For both generating sets one can restrict to root vectors 
$\xi$ and $\eta$.

%\bre{gen_set} Consider the subset $A_1$ of $R_q[U^w_+]$ 
%consisting of $c^{w, \om_i}_\eta$ where $i = 1, \ldots, r$ 
%and $\eta$ is a root vector of $V_w(\om_i)^*$. It is clear that 
%the generating set for the ideal $J_w(y)$ from \thref{main4} 
%is contained in $A_1 \cap J_w(y)$. In fact, these two sets 
%are different in general. To see this, recall \cite{G} 
%that for general $y \leq w$ the ideals $Q(y)^- + Q(w)^+$ 
%of $R^+$ are not prime (cf. \S 3.1 for the definitions 
%of $Q(y)^-$ and $Q(w)^+$). Only the saturated ideals 
%\[
%Q(y,w) = \{ c \in R^+ \mid \exists \la \in P_+
%\; \mbox{such that} \; c_w^\la a \in Q(y)^- + Q(w)^+ \}
%\] 
%are prime \cite[Proposition 6.8]{G}.
%Using $J_w(y) = \varphi_w(Q(y)_w^- + Q_w(w)^+)$,
%it is easy to see that 
%\[
%A_1 \cap J_w(y) = 
%\{ c^{w, \om_i}_{ \xi|_{V_w(\om_i)} } \mid 
%\xi \; \mbox{is a root vector of} \; V(\om_i)^* \; 
%\mbox{and} \; 
%c^{\om_i}_\xi \in Q(y,w) \}.
%\]
%In other words $A_1 \cap J_w(y)$ consists of those
%$c^{w, \om_i}_{ \xi|_{V_w(\om_i)} }$ where $\xi$ is a root vector of $V(\om_i)^*$
%for which there exists $\la \in P_+$ such that
%$\xi \otimes \xi_{w, \la} \in (V_w(\om_i + \la) 
%\cap \UU_- T_y v_{\om_i + \la})^\perp$, cf. \S \ref{2.5}.
%
%This set contains the generating set from \thref{main4}, but is
%in general bigger. Example 6.9 of Gorelik \cite{G} shows that this is 
%the case already for $\g = {\mathfrak{sl}}_3$ and 
%the ideal $I_{s_1 s_2}(s_1)$ of $R_q[U^{s_1 s_2}_+]$.
%\ere
%%%%%%%%%%%%%%%%%%%%%%%%%%%%%%%%%%%%%%%%%%%%%%%%%%%%%%%%%%%%%%%
\sectionnew{Results for the underlying Poisson structures}
\label{semi}
In this section we prove results for the underlying Poisson 
structures for the algebras $\UU^w_-$ which
are Poisson analogs of Theorems \ref{tmain2}, \ref{tmain3}
and \ref{tmain4}.
%%%%%%%%%%%%%%%%%%%%%%%%%%%%%%%%%%%%%%%%
\subsection{}
\label{4.1} Let $G$ be a simply connected complex semisimple Lie group and 
$\g = \Lie G$. Fix a pair of opposite Borel subgroups
$B_\pm$. Let $T= B_+ \cap B_-$ be the corresponding 
maximal torus of $G$, and $U_\pm$ be unipotent radicals 
of $B_\pm$. Let $W$ be the Weyl group of $G$. For all $w \in W$
fix representatives $\dot{w}$ in $N(T)/T$.  
Here $N(T)$ denotes the normalizer of $T$ in $G$.  

Denote by $\De_+$ the set of positive roots of $\g$. 
Fix root vectors $x^\pm_\al \in \g^{\pm \al}$, $\al \in \De_+$,
normalized by
\[
\lcor x^+_\al , x^-_\al \rcor = 1
\]
where $\lcor.,.\rcor$ denotes the Killing form on $\g$. 
Define the bivector field 
\[
\pi = - \sum_{\al \in \De_+} \chi(x^+_\al) \wedge \chi(x^-_\al) 
\]
on the flag variety $G/B_+$. Here $\chi \colon \g \to Vect(G/B_+)$
refers to the infinitesimal action of $\g$ on $G/B_+$. It is well 
known that $\pi$ is a Poisson structure on $G/B_+$, 
see e.g. \cite{GY} for details. The group $T$ acts on $(G/B_+, \pi)$
by Poisson automorphisms.

For $y_-, y_+ \in W$ define
\[
R_{y_-, y_+} = B_- y_- \cdot B_+ \cap B_+ y_+ \cdot B_+ \subset G/B_+.
\]
This intersection is nontrivial if and only if $y_- \leq y_+$
in which case it is irreducible \cite{De}. 

The following Proposition follows from \cite[Theorem 4.14]{EL}
of Evens and Lu, and \cite[Theorem 0.4]{GY} of Goodearl and 
the author.

\bpr{GY} The $T$-orbits of symplectic leaves of $(G/B_+, \pi)$ 
are precisely 
the intersections $R_{y_-, y_+}$, for $y_\pm \in W$, $y_- \leq y_+$.
\epr  

The closure relation between symplectic leaves is described by 
the well known fact that:
\[
\ol{R_{y_-, y_+}} = \bigsqcup \{ R_{y'_-, y'_+} \mid y'_\pm \in W, 
y_- \leq y'_- \leq y'_+ \leq y_+ \}.
\] 
%%%%%%%%%%%%%%%%%%%%%%
\subsection{}
\label{4.2}
From now on we fix an element $w \in W$. The Schubert (Bruhat) cell translate 
$w B_- \cdot B_+ \subset G/B_+$ is an open subset of $G/B_+$ 
and is thus a Poisson variety with the restriction of $\pi$. 

\ble{int} For $y_\pm \in W$ the intersection 
$w B_- \cdot B_+ \cap R_{y_-, y_+} \subset G/B_+$ is nonempty if and only if 
$y_- \leq w \leq y_+$. In the case when it is nontrivial, 
it is a dense subset of $R_{y_-, y_+}$. 
\ele
\begin{proof} The second statement holds because $w B_- \cdot B_+$ 
is a Zariski open subset of $G/B_+$ and $R_{y_-, y_+}$ are all irreducible.

Assume that $w B_- \cdot B_+ \cap R_{y_-, y_+}$ is nonempty. Then
\begin{multline}
w B_- B_+ \cap B_- y_- B_+ \neq \emptyset \Rightarrow
w B_- \cap B_- y_- B_+ \neq \emptyset \Rightarrow \\
B_- w B_- \cap B_- y_- B_+ \neq \emptyset \Rightarrow w \geq y_-.
\end{multline}
Analogously $w B_- B_+ \cap B_+ y_+ B_+ \neq \emptyset$ implies
$w \leq y_+$. Therefore, if 
$w B_- \cdot B_+ \cap R_{y_-, y_+}$ is nonempty, 
then $y_- \leq w \leq y_+$.

Now assume that that $y_- \leq w \leq y_+$. Analogously
we get $wB_- B_+ \cap B_+ y_+ B_+ \neq \emptyset$
and $w B_- B_+ \cap B_- y_- B_+ \neq \emptyset$.
Let $w b_- \in wB_- B_+ \cap B_+ y_+ B_+$
for some $b_- \in B_-$.
Since $wB_- B_+ \cap B_+ y_+ B_+$ is 
invariant under the left action of 
$B_+ \cap w B_- w^{-1}$
\begin{equation}
\label{incl1}
wB_- B_+ \cap B_+ y_+ B_+ \supset w (B_- \cap w^{-1} B_+ w) b_-.
\end{equation} 
Because $w B_- B_+ \cap B_- y_- B_+ \neq \emptyset$, 
$w B_- \cap B_- y_- B_+ \neq \emptyset$.
The latter set is invariant under the left action of 
$B_- \cap w B_- w^{-1}$ and 
$B_- = (B_- \cap w^{-1} B_- w)(B_- \cap w^{-1} B_+ w)$, thus
$w B_- \cap B_- y_- B_+ \neq \emptyset$ implies
\begin{equation}
\label{incl2}
\left( w (B_- \cap w^{-1} B_+ w) b_- \right) \cap B_- y_- B_+
\neq \emptyset.
\end{equation}
Then \eqref{incl1} and \eqref{incl2} imply
\[
B_- y_- B_+ \cap wB_- B_+ \cap B_+ y_+ B_+ \neq \emptyset.
\]
\end{proof}

\bpr{leaves} (1) The $T$-orbits of symplectic leaves of 
$(w B_- \cdot B_+, \pi)$ are the intersections
\[
S_w(y_-, y_+) = w B_- \cdot B_+ \cap R_{y_-, y_+} \quad
\mbox{for} \; (y_1, y_2) \in W \wDia W. 
\]
Their Zariski closures are given by
\[
\ol{S_w(y_-, y_+)} = \bigsqcup \{ S_w(y'_-, y'_+) \mid y'_\pm \in W, 
y_- \leq y'_- \leq w \leq y'_+ \leq y_+ \}.
\]

(2) Let $(y_-, y_+) \in W \wDia W$. For each symplectic leaf $\SS$ 
of $R_{y_-, y_+}$ the intersection $\SS \cap w B_- \cdot B_+$ is nontrivial
and is a symplectic leaf of $S_w(y_-, y_+)$. All symplectic leaves of 
$S_w(y_-, y_+)$ are obtained in this way. 
\epr
\begin{proof} Let $\SS$ be a symplectic leaf of $R_{y_-, y_+}$ 
and $(y_-, y_+) \in W \wDia W$. \prref{GY} implies that  
$R_{y_-, y_+} = T \cdot \SS$. Since the intersection $w B_- \cdot B_+ \cap R_{y_-, y_+}$
is nonempty there exists $t \in T$ such that $w B_- \cdot B_+ \cap t \SS \neq \emptyset$. 
But $w B_- \cap B_+$ is $T$-stable, so $w B_- \cdot B_+ \cap \SS \neq \emptyset$.
The complement of $w B_- \cdot B_+ \cap \SS$ in $\SS$ has real codimension at least 2.
Thus $w B_- \cdot B_+ \cap \SS$ is connected and is a symplectic leaf of 
$w B_- \cdot B_+$. Obviously all symplectic leaves of $w B_- \cdot B_+$ are 
obtained in this way, which completes the proof of (2). It is clear that 
\[
S_w(y_-, y_+) = T \cdot (w B_- \cdot B_+ \cap \SS).
\]
This implies (1).
\end{proof} 
%%%%%%%%%%%%%%%%%%%%%%%%%%%%%
\subsection{}
\label{4.3}
Denote
\begin{equation}
\label{Uw}
U^w_+ = U_+ \cap w U_- w^{-1} \; \; 
\mbox{and} \; \; \n^w_+ = \Lie U^w_+.
\end{equation}
Identify 
\begin{equation}
\label{iw} 
i_w \colon U^w_+ \cong B_+ w \cdot B_+ \subset G/B_+, \quad 
i_w(u) = uw \cdot B_+, \; 
u \in U^w_+.
\end{equation}  
Observe that 
$B_+ w \cdot B_+ = U^w_+ w \cdot B_+ = w(w^{-1} U^w_+ w) \cdot B_+$ lies 
inside $w B_- \cdot B_+$. \prref{GY} implies that 
$B_+ w \cdot B_+$ is a complete Poisson (locally closed) subset 
of $G/B_+$. Denote the Poisson structure 
\[
\pi_w = i_w^{-1}( \pi|_{B_+ w \cdot B_+})
\]
on $U^w_+$.

Consider the conjugation action of $T$ on $U^w_+$. It preserves $\pi_w$
since $i_w$ intertwines it with the canonical left action of $T$ on $G/B_+$.  

\bco{Uwleaves} The $T$-orbits of symplectic leaves
of $(U^w_+, \pi_w)$ are parametrized by $y \in W^{\leq w}$:
\begin{equation}
\label{mapw}
y \in W^{\leq w} \mapsto S_w(y) : = i_w^{-1}(R_{y, w})=
U^w_+ \cap B_- y B_+ w^{-1}.  
\end{equation} 
Moreover 
\[
\ol{S_w(y)} = \bigsqcup \{ S_w(y') \mid y' \in W^{\leq w}, y' \geq y \}.
\]  
In particular, \eqref{mapw} is an isomorphism of posets from $W^{\leq w}$
with the inverse Bruhat order to the underlying poset of the 
stratification of $(U^w_+, \pi_w)$ into $T$-orbits of symplectic 
leaves.
\eco
The Corollary follows from \prref{GY} since $i_w$ is $T$-equivariant and 
$B_+ w \cdot B_+$ is a complete Poisson (locally closed) subset 
of $w B_- \cdot B_+$.
%%%%%%%%%%%%%%%%%
\subsection{}
\label{4.4} 
The irreducible finite dimensional representations of $G$ are
parametrized by its set of positive dominant weights $P_+$.
Denote by $L(\la)$ the corresponding $G$-module. Let 
$d^\la_{\zeta, u} \in \Cset[G]$ be the 
the matrix coefficient corresponding to 
$\zeta \in L^*(\la)$ and $u \in L(\la)$. Then
\[
\Cset[G] = \Span \{ d^\la_{\zeta, u} \mid 
\la \in P_+, u \in L(\la),
\zeta \in L(\la)^* \}.
\]
We will denote the root spaces of a $G$-module $M$ by $M_\mu$.
For each $\la \in P_+$ fix a highest weight vector 
$u_\la$ of $L(\la)$ and a dual vector 
$\zeta_\la \in L^*(\la)_{-\la}$ normalized
by $\lcor \zeta_\la, u_\la \rcor =1$. 
Denote 
\[
d^\la_w = d^\la_{\dot{w} \zeta_\la, u_\la} 
\; \; \mbox{and} \; \;
d_w = \{ d^\la_w \mid \la \in P_+ \}.
\]
Then for $w B_- B_+ \subset G$
\[
\Cset[wB_-B_+] = \Cset[G] [d_w^{-1}].
\]
Identify 
\begin{equation}
\label{inv}
\Cset[w B_- \cdot B_+] \cong \Cset[w B_- B_+]^{B_+},
\end{equation}
where $(.)^{B_+}$ refers to the ring of invariant 
functions with respect to the right action of $B_+$ on $G$.
One verifies that under the isomorphism \eqref{inv}
\begin{equation}
\label{ringwB}
\Cset[w B_- \cdot B_+] = \{
d^\la_{\zeta, u_\la}/d^\la_w  \mid 
\la \in P_+, \zeta \in L(\la)^* \}.
\end{equation}
Analogously to \eqref{Rw0} one does not need to take
span in \eqref{ringwB}.
%%%%%%%%%%%%%%%%%%%%%%%
\subsection{}
\label{4.5} 
Denote $\n_\pm = \Lie U_\pm$. For $y \in W$, define the ideals 
\[
\wt{Q}(y)_w^\pm = \{d^\la_{\zeta, u_\la}/d^\la_w
\mid 
\la \in P_+, \xi \in (\UU(\n_\pm) y v_\la)^\perp
\subset L(\la)^* \}
\]
of $\Cset[wB_- \cdot B_+]$.

\bpr{ideal1} The vanishing ideal of the Zariski closure 
of $S_w(y, w)$ in $w B_- \cdot B_+$ is 
\begin{multline*}
\VV \left( Cl_{w B_- \cdot B_+}(S_w(y, w)) \right) = 
\wt{Q}(y)_w^- + \wt{Q}(w)_w^+
\\
= \{ d^\la_{\zeta, u_\la} / d^\la_w
\mid
\la \in P_+, 
\zeta \in ( \UU( \n_- ) y v_\la 
\cap \UU( \n_+ ) w v_\la)^\perp \subset L(\la)^*
\}.
\end{multline*}
\epr

\begin{proof} 
The ideal $\wt{Q}(y)_w^-$ is the vanishing ideal 
of $Cl_{w B_- \cdot B_+} (w B_- \cdot B_+ \cap B_- y \cdot B_+)$, 
in particular it is prime. Indeed 
\[
d_{\zeta, u_\la}^\la /d_w^{\la} \in  
\VV \left( Cl_{w B_- \cdot B_+}( w B_- \cdot B_+ \cap B_- y \cdot B_+)  
\right)
\]
if and only if $\lcor \zeta, (w B_- B_+ \cap B_- y B_+ u_\la) \rcor$
which is equivalent to $\lcor \zeta, B_- y B_+ u_\la \rcor$
and to $\zeta \in (\ol{\UU}_- u_\la)^\perp$ since
$w B_- B_+ \cap B_- y B_+$ is dense in $B_- y B_+$. Analogously,
in \S \ref{4.6} we verify that $\wt{Q}(w)_w^+$ is the 
vanishing ideal of $B_+ w \cdot B_+$ in 
$\Cset[w B_- \cdot B_+]$. 

Ramanathan proved \cite[Corollary 1.10 and Theorem 3.5]{R} 
that the scheme theoretic intersection 
of the opposite Schubert varieties $\ol{B_+ w \cdot B_+}$ and
$\ol{B_- y \cdot B_+}$ in $G/B_+$ is reduced. Therefore the same 
is true for the scheme theoretic intersection of 
$B_+ w \cdot B_+ = w B_- \cdot B_+ \cap \ol{B_+ w \cdot B_+}$ and 
$Cl_{w B_- \cdot B_+} (w B_- \cdot B_+ \cap B_- y \cdot B_+)=$
$w B_- \cdot B_+ \cap \ol{B_- y \cdot B_+}$.
This implies the statement of the Proposition.  
\end{proof}
%%%%%%%%%%%%%%%%%%%%%%%%
\subsection{}
\label{4.6}
For $\la \in P_+$ consider the $\UU(\b_+)$ submodules 
$L_w(\la) = \UU(\b_+) \dot{w} v_\la = \UU(\n^w_+) \dot{w} v_\la$ 
of $L(\la)$ (cf. \eqref{Uw}) called Demazure modules, 
where $\b_\pm = \Lie B_\pm$.
  
Each $\eta \in L_w(\la)^*$ gives rise to a regular function 
$d^{w, \la}_{\eta}$ on $U^w_+$, 
$d^{w, \la}_{\eta}(u) = \lcor \eta, u \dot{w} \zeta_\la \rcor$,
$u \in U^w_+$. One has
\begin{equation}
\label{Demfunct}
\Cset[U^w_+] =
\{ d^{w, \la}_{\eta} \mid \la \in P_+, \eta \in L_w(\la)^* \}.
\end{equation}
Let us trace back \eqref{Demfunct} to \eqref{ringwB}.
The composition of the isomorphism $i_w \colon U^w_+ \cong B_+ w \cdot B_+$
and the embedding $B_+ w \cdot B_+ \hookrightarrow w B_- \cdot B_+$ 
give rise to the embedding
\begin{equation}
\label{jw}
j_w \colon U^w_+ \hookrightarrow w B_- \cdot B_+, \quad
j_w(u) = u w B_+, u \in U^w_+.
\end{equation}
In terms of \eqref{ringwB} and \eqref{Demfunct} $j^*_w$ is given by
\begin{equation}
\label{jw*}
j^*_w(   d^\la_{\zeta, u_\la} / d^\la_w  ) = 
d^{w, \la}_{ \zeta|_{L_w(\la)}  }, \quad \la \in P_+, 
\zeta \in L(\la)^*.
\end{equation}
In particular, the kernel of $j^*_w$ is
\begin{multline}
\label{kerj}
\ker j^*_w = \VV(B_+ w \cdot B_+) = \VV( \ol{S_w(1, w)} ) \\
=\{d^\la_{\zeta, u_\la}/d^\la_w
\mid \la \in P_+, \zeta \in (\UU(\n_+) w v_\la)^\perp 
\subset L(\la)^* \},
\end{multline}
cf. \prref{leaves} and \prref{ideal1}.
\bth{mainP1} For all $y \in W^{\leq w}$ the vanishing ideal of 
the Zariski closure of the symplectic leaf $S_w(y)$ in 
$(U^w_+, \pi_w)$ is 
\[
\VV(\ol{S_w(y)}) =  \{ d^{w,\la}_\eta \mid
\eta \in (L_w(\la) \cap \UU(\n_-) y u_\la)^\perp 
\subset L_w(\la)^* \}.
\]
\eth
\begin{proof} Clearly $j_w( \ol{S_w(y)} )= \ol{ S_w(y,w) }$, 
cf. \prref{leaves} and \coref{Uwleaves}
(equivalently one can use that $j_w$ is closed). Thus
$\VV(\ol{S_w(y)}) = j^*_w (\VV( \ol{S_w(y,w)} ))$ and the Theorem
follows from \prref{ideal1} and \eqref{jw*}.
\end{proof}

We complete this subsection with a proof of the fact that 
for $y \in W^{\leq w}$ the ideal 
$\wt{Q}(y, w)_w= \wt{Q}(y)_w^- + \wt{Q}(w)_w^+$
is prime. 
%%%%%%%%%%%%%%%%%%%%%%%%%%%%%%%%%%%%%%
\subsection{}
\label{4.7}
The following Theorem is a Poisson analog of \thref{main4}. 

\bth{mainP3} For all $y \in W^{\leq w}$ the vanishing ideal of 
the Zariski closure of the $T$-orbit of symplectic leaves 
$S_w(y)$ in $(U^w_+, \pi_w)$ is generated by $d^{w,\om_i}_\eta$
where $i = 1, \ldots, r$ and 
$\eta \in (L_w(\om_i) \cap \UU(\n_-) y u_{\om_i})^\perp \subset L_w(\om_i)^*$
is a root vector.
\eth
\begin{proof} Consider the algebra $\Cset[G]^{U_+}$ of right 
$U^+$-invariant functions on $G$. It is spanned by the matrix 
coefficients $d^\la_{\zeta, u_\la}$, $\la \in P_+$, $\zeta \in L(\la)^*$.
Kempf and Ramanathan proved \cite[Theorem 3(i)]{KR} that Schubert varieties 
are linearly defined. This implies that for all $y \in W$
\[
\Span \{ d^\la_{\zeta, u_\la} \mid
\la \in P_+, \zeta \in (\UU(\n_-) \dot{y} u_\la)^\perp \}
= \sum_{i=1}^r \{ d^\la_{\zeta, u_{\om_i}} \mid
\zeta \in (\UU(\n_-) \dot{y} u_{\om_i})^\perp \} \Cset[G]^{U_+}.
\]
Then inside $\Cset[wB_- \cdot B_+]$ one has 
\begin{multline*}
\Span \{ d^\la_{\zeta, u_\la}/d^\la_w \mid
\la \in P_+, \zeta \in (\UU(\n_-) \dot{y} u_\la)^\perp \}
\\ = \sum_{i=1}^r \{ d^\la_{\zeta, u_{\om_i}}/d^{\om_1}_w 
\mid \zeta \in (\UU(\n_-) \dot{y} u_{\om_i})^\perp \} 
\Cset[wB_- \cdot B_+],
\end{multline*}
recall \eqref{ringwB}. Now the Theorem follows from
$\VV(\ol{S_w(y)}) = j^*_w (\VV( \ol{S_w(y,w)} ))$
and \prref{ideal1}.
\end{proof}
%%%%%%%%%%%%%%%%%%%%%%%%%%%%%%%%%%%%%%%%%%%%%%%%%%%%%%%%%%%%%%%
\sectionnew{Quantum matrices}
\label{RqM}
%%%%%%%
\subsection{}
\label{5.1}
Throughout this section we fix two positive integers
$m$ and $n$. Let $G = SL_{m+n}(\Cset)$ and $B_\pm$ 
be its standard Borel subgroups.

Denote by $w^\ci_{m+n}$ the longest element 
of $S_{m+n}$. For each $k \leq m+n$ denote by 
$w^\ci_k$ and $w^{\ci r}_k$ the longest
elements of $S(\{1, \ldots, k \}) \subseteq S_{m+n}$
and $S(\{m+n-k+1, \ldots, m+n \}) \subseteq S_{m+n}$, 
respectively.  

Denote the Coxeter element $c = (1 2 \ldots m+n) \in S_{m+n}$.
Then
\begin{equation}
\label{Coxeter}
c^m = w^\ci_m w^{\ci r}_n w^\ci_{m+n}.
\end{equation}

In \S 5.1-5.3 we will apply the results of the previous 
Section to the case $G = SL_{m+n}(\Cset)$, 
$\g = {\mathfrak{sl}}_{m+n}(\Cset)$ and $w = c^m$. 
All notation  $L(\om_k)$, $L_w(\om_k)$, $\n_+$, 
$U^w_+$, $\pi_w$ will refer to this case. 

For two integers $k \leq l$ set $\ol{k,l}=\{k, \ldots, l\}$.   
%%%%%%%%%%%%%%%%%%%%%%%%%%%%%%%%%%
\subsection{}
\label{5.2}
The matrix affine Poisson space is the complex affine space $\Mmn$
consisting of rectangular matrices of size $m \times n$
equipped with the quadratic Poisson structure
\begin{equation}
\lb{pi_mn}
\pi_{m,n} = \sum_{i, k=1}^m \sum_{j, l=1}^n
(\sign(k-i) + \sign(l-j)) x_{il} x_{kj}
\frac{\partial}{\partial x_{ij}} \wedge
\frac{\partial}{\partial x_{kl}},
\end{equation}
where $x_{ij}$ are the standard coordinate functions on $\Mmn$.

One has, cf. \cite[Proposition 3.4]{BGY}, \cite[(3.11)]{GSV}, 
\cite[Proposition 1.6]{GY}: 

\bpr{Matisom}
The map $f \colon (\Mmn, \pi_{m,n}) \to (U^{c^m}_+, \pi_{c^m})$
given by 
\[
f(x) = \left( \begin{smallmatrix} I_m & w_m^\ci x \\ 0 & I_n
\end{smallmatrix} \right)
\]
is an isomorphism of Poisson varieties, where 
$U^{c^m}_+ \subset SL_{m+n}(\Cset)$ is given by \eqref{Uw} 
\epr

Here, for $w \in S_m$ we denote by the same letter 
the corresponding permutation matrix in $GL_m(\Cset)$. 

Define the torus $T := \Cset^{m+n-1}$ and view it 
as pairs of diagonal matrices $(A, B)$ of size $m \times m$ and 
$n \times n$ with $\det (A) \det (B) =1$. It acts on $\Mmn$ by 
$(A, B) \cdot X = A X B^{-1}$, $X \in \Mmn$. 
The Poisson structure $\pi_{m,n}$ is invariant under 
the action of $T$ and $f$ intertwines it with the conjugation 
action of the standard torus of $SL_{m+n}(\Cset)$ 
on $U^{c^m}_+$, see \S 4.3. For $y \in S_{m+n}^{ \leq c^m}$ 
denote
\[
S(y) = f^{-1}(B_- y B_+ c^{-m} ),
\]
where $B_\pm$ refer to the standard Borel subgroups of 
$SL_{m+n}(\Cset)$.

\bco{Mleaves} \cite[Theorem A]{BGY} The $T$-orbits of 
symplectic leaves of $(\Mmn, \pi_{m,n})$ are 
$S(y)$, $y \in S_{m+n}^{ \leq c^m}$. Their Zariski 
closures are given by 
\[
\ol{S(y)} = \bigsqcup_{ y \leq y' \leq c^m} S(y'). 
\]
\eco
%%%%%%%%%%%%%%%%%%%%%%%%%%%%%%%%%%
\subsection{}
\label{5.3}
Denote by $L$ the vector representation of $SL_{m+n}(\Cset)$ 
with standard basis $\{ u_1, \ldots, u_{m+n} \}$ (such that 
$E_{ij} u_q = \delta_{jq} u_i$). The fundamental 
representations of $SL_{m+n}(\Cset)$ are 
$L(\om_k) \cong \wedge^k L$, $k=1, \ldots, m+n-1$. 
They have bases 
\[
u_I = u_{i_1} \wedge \ldots \wedge u_{i_k}, \quad 
I = \{ i_1 < \cdots < i_k \} \subset \ol{1,m+n}.
\]

For a subset $I \subseteq \ol{1,m+n}$ denote
\begin{equation}
\label{p1}
p_1 (I) = I \cap \ol{1,m} \; \; \mbox{and} \; \; 
p_2 (I) = I \cap \ol{m+1, m+n}.
\end{equation}

Consider the partial order on $\{ I \subseteq \ol{1,m} \mid |I|= k\}$:
for $I= \{ i_1 < \cdots < i_k \}, J= \{ j_1 < \cdots < j_k \} 
\subseteq \ol{1,m+n}$, $I \leq J$ if 
$i_l \leq j_l$ for all $l=1, \ldots, k$.

For $J_1 \subseteq \ol{1,m}$, $J_2 \subseteq \ol{1,n}$, 
$|J_1|=|J_2|$ 
denote by $\De_{J_1,J_2}(x)$ the corresponding minor of 
$x \in \Mmn$.

Let $I \subset \ol{1,m+n}$. If $k \in \ol{1,n}$, then 
$I \leq c^m(\ol{1,k}) = \ol{m+1, m+k}$ implies 
$p_2(I) \subseteq \ol{m+1, m+k}$. If $k \in \ol{n+1,m+n-1}$, 
then $I \leq c^m(\ol{1,k}) = \ol{1,k-n} \sqcup \ol{ m+1, m+n}$ 
implies $p_1(I) \supseteq \ol{1, k-n}$. For 
$y \in S_{m+n}^{\leq c^m}$ let $\AA(y)$ be the union of the sets of minors
\begin{equation}
\label{det1}
\Delta_{w_m^\ci(p_1(I)), (\ol{m+1, m+k} \backslash p_2(I)) -m}
\end{equation}
for $k \in \ol{1,n}$, $I \subseteq \ol{1, m+n}$,
$|I|=k$, $I \leq c^m( \ol{1,k})$, 
$I \ngeq y(\ol{1,k})$
and
\begin{equation}
\label{det2}
\Delta_{w_m^\ci(p_1(I) \backslash \ol{1,k-n}), 
(\ol{m+1, m+n} \backslash p_2(I)) -m}
\end{equation}
for
$k \in \ol{n+1,m+n-1}$, $I \subset \ol{1, m+n}$,
$|I|=k$, $I \leq c^m( \ol{1,k})$, $I \ngeq y(\ol{1,k})$.
In \eqref{det1}--\eqref{det2}
$-m$ means subtracting $m$ from each element 
of the set. Both sets of minors \eqref{det1}--\eqref{det2} 
can be uniformly 
described by the less explicit formula \\
$\Delta_{p_1(I) \backslash p_1(c^m(\ol{1,k})) , 
(p_2 (c^m(\ol{1, k})) \backslash p_2(I)) -m}$,
$k \in \ol{1,m+n-1}$.

\bth{IdMat} For all $y \in S_{m+n}^{\leq c^m}$ the vanishing ideal 
of the Zariski closure of the $T$-orbit of symplectic leaves $S(y)$ in 
$(\Mmn, \pi_{m,n})$ is generated by the minors 
in $\AA(y) \subset \Cset[\Mmn]$. 
\eth

Functions cutting the closures $\ol{S(y)}$ were previously 
obtained by Brown, Goodearl and the author in 
\cite[Theorem 4.2]{BGY}. Goodearl, Launois and Lenagan
\cite{GLL1} independently find all minors that belong 
to the vanishing ideal of the Zariski closure of 
any $T$-orbit of symplectic leaves $S(y)$ in 
$(\Mmn, \pi_{m,n})$. 

\begin{proof} Observe that the Demazure 
module $L_{c^m}(\om_k)$ is given by 
\[
L_{c^m}(\om_k) = \Span \{ u_I \mid I \subset \ol{1, m+n}, |I|=k, 
I \leq c^m(\ol{1,k}) \}.
\]
Denote by $\{ \zeta_I \mid I \subset \ol{1, m+n}, |I|=k \}
\subset L(\om_k)^*$ the dual basis to $\{ u_I \}$.
We can identify
\begin{equation}
\label{Demident}
L_{c^m}(\om_k)^* \cong \Span \{ \zeta_I \mid I \subset \ol{1, m+n}, |I|=k, 
I \leq c^m(\ol{1,k}) \} \subset L(\om_k)^*.
\end{equation}
Then 
\[
L_{c^m}(\om_k) \cap \UU(\n_-) \dot{y} u_{\ol{1,k}} 
= \Span \{ u_I \mid I \subseteq \ol{1, m+n}, |I|=k, 
y(\ol{1,k}) \leq I \leq c^m(\ol{1,k})\},
\]
where $\n_-$ denotes the nilpotent subalgebra of ${\mathfrak{sl}}_{m+n}$
consisting of lower triangular matrices.
Under the identification \eqref{Demident} the orthogonal complement 
$(L_{c^m}(\om_k) \cap \wt{\UU}_- \dot{y} u_{\ol{1,k}})^\perp$ in 
$L_{c^m}(\om_k)^*$ is 
\[
(L_{c^m}(\om_k) \cap \wt{\UU}_- \dot{y} u_{\ol{1,k}})^\perp
= \Span \{ \zeta_I \mid 
I \subseteq \ol{1, m+n}, |I|=k, 
I \leq c^m(\ol{1,k}), I \ngeq y(\ol{1,k}) \}.
\]
\thref{mainP3} implies that the vanishing ideal of the 
Zariski closure of $S(y)$ in $\Mmn$ is generated 
by $d_{\zeta_I}^{c^m, \om_k}(f(x) )$, 
$k \in \ol{1, m+n -1}$, $I \subset \ol{1, m+n}$, 
$|I|=k$, $I \leq c^m(\ol{1,k}), I \ngeq y(\ol{1,k})$.
It is straighforward to check that 
\[
d_{\zeta_I}^{c^m, \om_k}(f(x) )
=
\begin{cases}
\Delta_{w_m^\ci(p_1(I)), (\ol{m+1, m+k} \backslash p_2(I)) -m}(x),  
& \mbox{if}\; 1 \leq k \leq n \\
\Delta_{w_m^\ci(p_1(I) \backslash \ol{1,k-n}), 
(\ol{m+1, m+n} \backslash p_2(I)) -m}(x),  
& \mbox{if}\; n+1 \leq k \leq m+n -1
\end{cases}
\]
This completes the proof of the Theorem.
\end{proof}
%%%%%%%%%%%%%%%%%%%%%%%%%%%%%%%%%%%
\subsection{}
\label{5.4}
The algebra of quantum matrices $R_q(\Mmn)$ is the 
$\KK$-algebra generated by $x_{ij}$, $i= 1, \ldots,m$,
$j = 1, \ldots, n$ subject to the relations
\begin{align*}
x_{ij} x_{lj} &= q x_{lj} x_{ij}, \quad \mbox{for} \; i < l, \\
x_{ij} x_{ik} &= q x_{ik} x_{ij}, \quad \mbox{for} \; j < k, \\ 
x_{ij} x_{lk} &= x_{lk} x_{ij}, \quad 
\mbox{for} \; i < l, j > k,\\
x_{ij} x_{lk} - x_{lk} x_{ij} &= (q-q^{-1}) x_{ik} x_{lj}, \quad 
\mbox{for} \; i < l, j<k,
\end{align*}
where $\KK$ is a field of characteristic 0 and $q \in \KK$
is transcendental over $\Qset$.
For $I = \{ i_1 < \cdots < i_k \} \subset \{1, \ldots, m\}$
and $J = \{ j_1 < \cdots < j_k \} \subset \{1, \ldots, n\}$
one defines the quantum minor $\De^q_{I,J} \in R_q(\Mmn)$
by 
\begin{align}
\nn
\De^q_{I,J} &= \sum_{\sig \in S_k} (-q)^{l(\sig)}
x_{i_1 j_{\sig(1)}} \ldots x_{i_k j_{\sig(k)}}
\\
&= \sum_{\sig \in S_k} (-q)^{-l(\sig)}
x_{i_k j_{\sig(k)}} \ldots x_{i_1 j_{\sig(1)}}.
\label{qminor}
\end{align}

The group $\Zset^{m+n}$ acts on $R_q[\Mmn]$ 
by algebra automorphisms by setting 
$(a_1, \ldots, a_m, b_1, \ldots, b_n) \cdot x_{ij} 
= q^{a_i - b_j} x_{ij}$ on the generators 
of $R_q[\Mmn]$.
%%%%%%%%%%%%%%%%%%%%%%%%%%%%%%%%%%
\subsection{}
\label{5.5}
In \S \ref{5.5}--\ref{5.7} we apply the results from 
Sect. \ref{hsp} to the particular case $\g = {\mathfrak{sl}}_{m+n}$, $w = c^m$. 
In particular $\UU_+$, $V(\om_k)$, $V_w(\om_k)$ refer to 
this situation. 

Consider the reduced decomposition
\begin{equation}
\label{w0}
w^\ci_{m+n} = s_1 (s_2 s_1) \ldots (s_{m+n-1} \ldots s_1).
\end{equation}
Denote the corresponding root vectors given by \eqref{rootv2} by
\[
Y_{1,2}; Y_{1,3}, Y_{2,3}; \ldots; Y_{1,m+n}, \ldots, Y_{m+n-1, m+n} 
\in \UU^{w^\ci_{m+n}}_+ = \UU_+
\]
and
\[
Y_{2,1}; Y_{3,1}, Y_{3,2}; \ldots; Y_{m+n,1}, \ldots, Y_{m+n, m+n-1} 
\in \UU^{w^\ci_{m+n}}_- = \UU_-
\]
in the plus and minus cases, respectively.
Then by \cite[Lemma 2.1.1]{MC} $Y_{i, i+1} = X_i^+$, $1 \leq i < m+n-1$ 
and for $i< j$ $Y_{ij}$ is recursively given by
\begin{equation}
\label{rec1}
Y_{ij} = Y_{i, j-1}Y_{j-1,j} - q^{-1} Y_{j-1,j} Y_{i, j-1}.
\end{equation}
Analogously one has that 
$Y_{i+1, i} = X_i^-$, $1 \leq i <m+n-1$ and for $j>i$ 
$Y_{ji}$ is recursively given by
\begin{equation}
\label{rec2}
Y_{ji} = Y_{j,j-1} Y_{j-1, i} - q Y_{j-1,i} Y_{j,j-1}.
\end{equation}
The expression
\[
c^m = (s_m \ldots s_1) (s_{m+1} \ldots s_2) \ldots (s_{m+n-1} \ldots s_n)
\]
is reduced since $l(c^m) = mn$. 
Denote the corresponding root vectors of $\UU^{c^m}_+$ by 
\[
X_{1,m+1}, \ldots, X_{m,m+1}; X_{1, m+2}, \ldots, X_{m,m+2}; \ldots; 
X_{1, m+n}, \ldots, X_{m,m+n}
\] 
and of $\UU^{c^m}_-$ by 
\[
X_{m+1,1}, \ldots, X_{m+1,m}; X_{m+2,1}, \ldots, X_{m+2,m}; \ldots; 
X_{m+n,1}, \ldots, X_{m+n,m}.
\] 

\ble{isom} (1) For all $i \in \ol{1,m}$, $j \in \ol{m+1,m+n}$ and
$i \in \ol{m+1, m+n}$, $j \in \ol{1,m}$:
\[
X_{ij} = T_{w^\ci_m} Y_{ij}.
\]

(2) The map $g \colon R_q[\Mmn] \to \UU^{c^m}_-$ given by
\[
x_{ij} \mapsto (-q)^{j+m-i-1}X_{j+m,i}, \quad i \in \ol{1,m}, j \in \ol{1,n}
\]
is an isomorphism of algebras.
\ele

\begin{proof}
The first part of (1) is \cite[Lemma 2.1.3 (3)]{MC}.
The second part of (1) is similar. M\'eriaux and 
Cauchon showed that $x_{ij} \mt Y_{i,j+m}$ defines an algebra isomorphism,
based on the Alev--Dumas result \cite{AD} that the Yamabe root vectors
of $\UU_q({\mathfrak{sl}}_{m+n})$ satisfy 
the relations for the standard generators of $R_q[\Mmn]$. 
Since $X_i^+ \mt X_i^-$
defines an isomorphism from $\UU_+$ to $\UU_-$ (such that 
$Y_{ij} \mt (-q)^{j-i-1} Y_{ji}$ for $i <j$) and 
$T_{w^\ci_m}$ is an (algebra) automorphism of $\UU_q(\g)$, the map $g$
is a homomorphism. It is an isomorphism because of the PBW basis 
part of \thref{DKP}. This proves (2).   
\end{proof}
%%%%%%%%%%%%%%%%%%%%%
\subsection{}
\label{5.6}
For $y \in S_{m+n}^{\leq c^m}$ let $\AA_q(y)$ be the union 
of the sets of quantum minors
\[
\Delta^q_{w_m^\ci(p_1(I)), (\ol{m+1, m+k} \backslash p_2(I)) -m}
\in R_q[\Mmn]
\]
for $k \in \ol{1,n}$, $I \subset \ol{1, m+n}$,
$|I|=k$, $I \leq c^m( \ol{1,k})$, 
$I \ngeq y(\ol{1,k})$
and
\[
\Delta^q_{w_m^\ci(p_1(I) \backslash \ol{1,k-n}), 
(\ol{m+1, m+n} \backslash p_2(I)) -m} \in R_q[\Mmn]
\]
for $k \in \ol{n+1,m+n-1}$, $I \subset \ol{1, m+n}$,
$|I|=k$, $I \leq c^m( \ol{1,k})$, $I \ngeq y(\ol{1,k})$, 
cf. \eqref{p1}. We refer the reader to 
\eqref{det1}--\eqref{det2} for a comparison to the 
Poisson case.

\bth{Qmatr} For all $y \in S_{m+n}^{\leq c^m}$ denote by
$I(y)$ the right ideal of $R_q[\Mmn]$ generated by $\AA_q(y)$. 

Then all ideals $I(y)$ are two sided, prime 
and $\Zset^{m+n}$-invariant. 
They exhaust all $\Zset^{m+n}$-primes 
of $R_q[\Mmn]$. The map $y \in S_{m+n}^{\leq c^m}
\mt I(y)$ is an isomorphism from the poset 
$S_{m+n}^{\leq c^m}$ to the poset of $\Zset^{m+n}$
invariant prime ideals 
of $R_q[\Mmn]$ ordered under inclusion.
\eth

\thref{Qmatr} is a corollary of Theorems \ref{tmain2}, 
\ref{tmain3} and \ref{tmain4} for the special case of 
the algebras $\UU^{c^m}_-$, cf. \S \ref{5.5}. 
Its proof will be given in \S \ref{5.7}.

The parametrization and poset structure of 
$\Zset^{m+n}$-primes is due to 
Launois \cite{La} who also proved that all
of them are generated by quantum minors. 
Our proof is independent. Generators 
for the $\Zset^{m+n}$-primes of $R_q[\Mmn]$ 
were only known in the case $m=n=3$ due 
to Goodearl and Lenagan \cite{GL}. Goodearl, Launois, and 
Lenagan have a recent independent approach 
constructing ideal generators in the
general case \cite{GLL}. 

Define the algebra
\[
\Lambda_q(\KK^{m+n}) = T_\KK(v_1, \ldots, v_{m+n}) / 
\lcor v_i v_j = - q^{-1} v_j v_i, i > j, 
v_i^2=0 \rcor
\]
where $T_\KK(.)$ refers to the tensor algebra over $\KK$.
It has a canonical structure of 
$\UU_q({\mathfrak{sl}}_{m+n})$-module algebra 
for the action:
\begin{align*}
Y_{ij} v_k &= \delta_{jk} v_i, 
\\
K_i v_k &= q^{a_{ik}} v_k, \quad
a_{ki} = 1 \; \mbox{if} \; k=i, \;
a_{ki}= -1 \; \mbox{if} \; k=i-1, 
\; a_{ki}= 0 \; \mbox{otherwise}.  
\end{align*}
Moreover $\Lambda_q(\KK^{m+n})$ is graded by $\deg v_i = 1$ and its
$k$-graded component is isomorphic to the 
fundamental representation $V(\om_k)$ of 
$\UU_q({\mathfrak{sl}}_{m+n})$ 
\begin{equation}
\label{isomqsl}
\Lambda_q(\KK^{m+n})_k \cong V(\om_k)
\end{equation}
for $k = 1, \ldots, m+n-1$. We will 
use the isomorphism \eqref{isomqsl} for the remainder of
this Section. Assuming it,
\[
v_I := v_{i_1} \ldots v_{i_k}, \quad
I=\{i_1 < \cdots< i_k\} \subset \ol{1,m+n}.
\]
is a basis of $V(\om_k)$.
Denote the dual basis of 
$V(\om_k)^*$ by $\{ \xi_I \}$. 
Since all root spaces of $V(\om_k)$ are one dimensional 
\eqref{Tcommut} implies
\begin{equation}
\label{Tv}
T_w v_I = b v_{w(I)}
\end{equation}
for some nonzero $b \in \KK$ (depending on $I$ and $w$). 
Recall the partial order on 
$\{ I \subset \ol{1,m+n} \mid |I|=k \}$ from 
\S \ref{5.3}. Then the Demazure module 
$V_w(\om_k)$ is given by 
\[
V_w(\om_k) = \UU_+ T_w v_{\ol{1,k}} = \Span
\{ v_I \mid I \subset \ol{1,m+n}, 
|I|=k, I \leq w(\ol{1,k}) \}.
\]  
Identify the dual space $V_w(\om_k)^*$ 
with 
\begin{equation}
\label{dualDemazure}
V_w(\om_k)^* \cong \Span
\{ \xi_I \mid I \subset \ol{1,m+n}, 
|I|=k, I \leq w(\ol{1,k}) \} 
\subset V(\om_k)^*.
\end{equation}
Under this identification the orthogonal
complement $(V_w(\om_k)\cap \UU_- T_y v_{\ol{1,k}})^\perp$
to $V_w(\om_k)\cap \UU_- T_y v_{\ol{1,k}}$ in $V_w(\om_k)^*$
is given by 
\[
(V_w(\om_k)\cap \UU_- T_y v_{\ol{1,k}})^\perp =
\Span \{ \xi_I \mid I \subset \ol{1,m+n}, |I|=k,
I \leq c^m(\ol{1,k}), I \ngeq y(\ol{1,k} \}
\]  
for all $y \in S_{m+n}^{\leq c^m}$.
%%%%%%%%%%%%%%%%%%%%%%%%%%%%%%%
\subsection{}
\label{5.7}
We have 
\begin{align*}
\RR^{c^m} = &\left( \exp_q(X_{m, m+n} \otimes X_{m+n,m}) \dots 
\exp_q(X_{1, m+n} \otimes X_{m+n,1}) \right) \dots
\\
&\left( \exp_q(X_{m, m+1} \otimes X_{m+1,m}) \dots 
\exp_q(X_{1, m+1} \otimes X_{m+1,1}) \right),
\end{align*}
recall \eqref{Rw}. From \eqref{Tcommut} one obtains that 
\begin{equation}
\label{Teq}
(T_{w_m^\ci}^{-1} (x) ). v = T_{w_m^\ci}^{-1} x T_{w_m^\ci} v
\end{equation}
for all $x \in \UU_q({\mathfrak{sl}}_{m+n})$, $v \in V(\om_k)$. 

Denote
\[
c_{I,J}^k = c^{\om_k}_{\xi_I, v_J}.
\]
Taking into account \leref{isom}, eqs. 
\eqref{Tv}, \eqref{Teq}, and the fact that $T_w$ is an 
algebra automorphism of $\UU_q({\mathfrak{sl}}_{m+n})$
for all $w \in W$, we obtain:
\begin{multline}
\label{identRR}
(c^k_{I,J} \otimes g^{-1}) (\RR^{c^m}) = b
\big(c^k_{w_m^\ci(I),w_m^\ci(J)} \otimes \id \big) \big[
\big( \exp_q(Y_{m, m+n} \otimes x_{m,n}) \dots \\ 
\exp_q( Y_{1, m+n} \otimes x_{1,n}) \big) \dots
\big( \exp_q(Y_{m, m+1} \otimes x_{m,1} ) \dots 
\exp_q(Y_{1, m+1} \otimes x_{1,1}) \big)
\big]
\end{multline}
for some nonzero $b \in \KK$.
\\
\noindent
{\em{Proof of \thref{Qmatr}.}} For $y \in S_{m+n}^{\leq c^m}$ 
define the right ideals
\begin{multline}
\label{Iwtide}
\wt{I}(y) =\big\{ (c^k_{I, c^m(\ol{1,k})} \otimes g^{-1} )(\RR^{c^m})
\, \big| \, k \in \ol{1, m+n-1}, I \subset \ol{1,m+n}, \\
|I|=k, I \leq c^m(\ol{1,k}), I \ngeq y(\ol{1,k}) \big\}
R_q[\Mmn].
\end{multline}

Theorems \ref{tmain2}, \ref{tmain3} and \ref{tmain4}
and \leref{isom} (2) imply that all ideals $I(y)$ 
are two-sided, prime and $\Zset^{m+n}$-invariant. 
Moreover they exhaust all $\Zset^{m+n}$-primes 
of $R_q[\Mmn]$ and the map $y \in S_{m+n}^{\leq c^m}
\mt I(y)$ is an isomorphism from the poset 
$S_{m+n}^{\leq c^m}$ to the poset of $\Zset^{m+n}$ primes 
of $R_q[\Mmn]$ ordered under inclusion.

We claim that for $k \in \ol{1,n}$, $I \subset \ol{1,m+n}$,
$|I|=k$, $I \leq c^m(\ol{1,k})$
\begin{equation}
\label{De1}
(c^k_{I, c^m(\ol{1,k})} \otimes g^{-1} )(\RR^{c^m}) =
b \Delta^q_{w_m^\ci(p_1(I)), (\ol{m+1, m+k} \backslash p_2(I)) -m}
\end{equation}
and 
for $k \in \ol{n+1,m+n-1}$, $I \subset \ol{1,m+n}$,
$|I|=k$, $I \leq c^m(\ol{1,k})$
\begin{equation}
\label{De2}
(c^k_{I, c^m(\ol{1,k})} \otimes g^{-1} )(\RR^{c^m}) = b
\Delta^q_{w_m^\ci(p_1(I) \backslash \ol{1,k-n}), 
(\ol{m+1, m+n} \backslash p_2(I)) -m}
\end{equation}
for some nonzero $b \in \KK$ depending on $k$ and $I$.
This implies that $I(y) = \wt{I}(y)$ for all 
$y \in S_{m+n}^{\leq c^m}$ and the statement of the Theorem.

Eqs. \eqref{De1} and \eqref{De2} are verified in a similar way.
We will restrict ourselves to \eqref{De1}. 
Let $i <j \in \ol{1,m+n}$. From \eqref{rec1} one checks 
inductively on $j-i$ that 
\begin{equation}
\label{Yid1}
Y_{ij} (u_{I'} u_j) = u_{I'} u_i
\end{equation}
for all $I' \subset \ol{1, j-1} \sqcup \ol{j+1, m+n}$, $|I'|=k-1$
and
\begin{equation}
\label{Yid2}
Y_{ij} (u_I) = 0
\end{equation}
for all $I \subset \ol{1, j-1} \sqcup \ol{j+1, m+n}$, $|I|=k$.

Now fix $k \in \ol{1,n}$ and $I \subset \ol{1,m+n}$ such that
$|I|=k$, $I \leq c^m(\ol{1,k})$.
Compute $w_m^\ci(I) = w_m^\ci(p_1(I)) \sqcup p_2(I)$ 
and $w_m^\ci c^m (\ol{1,k}) = \ol{m+1,m+k}$. Then 
\[
w_m^\ci(I) \cap w_m^\ci c^m (\ol{1,k}) = \ol{m+1, m+k} \cap p_2(I)
\]
and 
\[
w_m^\ci(I) \backslash w_m^\ci c^m (\ol{1,k}) = w_m^\ci(p_1(I)),
\quad
w_m^\ci c^m (\ol{1,k}) \backslash w_m^\ci(I) = 
\ol{m+1, m+k} \backslash p_2(I).
\]
Denote
\[
w_m^\ci(p_1(I)) = \{ i_1 < \ldots < i_l \}, \quad
\ol{m+1, m+k} \backslash p_2(I) = \{ j_1 +m < \ldots < j_l + m \}.
\]
Eqs. \eqref{identRR} and \eqref{Yid1}--\eqref{Yid2} imply
\begin{multline*}
(c^k_{I, c^m(\ol{1,k}) } \otimes g^{-1} )(\RR^{c^m}) = 
\\
b_1
\sum_{\sigma \in S_l} 
x_{i_{\sigma(l)}, j_l} \ldots x_{i_{\sigma(1)}, j_1}
\lcor \xi_{i_1,\ldots, i_l}, 
Y_{i_{\sigma(l)}, j_l+m} \ldots Y_{i_{\sigma(1)}, j_1+m} 
v_{j_l+m} \ldots v_{j_1+m} \rcor
\end{multline*}
for some nonzero $b_1 \in \KK$. Using the fact that 
$v_i v_j = -q^{-1} v_j v_i$ for $i >j$, eq. \eqref{Yid1} and the 
fact that for a permutation $\sigma$, $l(\sigma)$ is equal to the 
number of its inversions we obtain
\begin{align*}
(c^k_{I, c^m(\ol{1,k}) } \otimes g^{-1} )(\RR^{c^m}) 
&= q^{-l(l-1)/2} b_1 \sum_{\sigma \in S_l} 
x_{i_{\sigma(l)}, j_l} \ldots x_{i_{\sigma(1)}, j_1}
\lcor 
\xi_{i_1,\ldots, i_l},
v_{i_{\sigma(1)}} \ldots v_{i_{\sigma(l)}} \rcor
\\
&= q^{-l(l-1)/2} b_1 \sum_{\sigma \in S_l} (-q)^{-l(\sigma)}
x_{i_{\sigma(l)}, j_l} \ldots x_{i_{\sigma(1)}, j_1}
\\
&= q^{-l(l-1)/2} b_1 
\Delta^q_{w_m^\ci(p_1(I)), (\ol{m+1, m+k} \backslash p_2(I)) -m}.
\end{align*}
This completes the proof of \eqref{De1} and the Theorem. 
\qed
%%%%%%%%%%%%%%%%%%%%%% References %%%%%%%%%%%%%%%%%%%%%%%%%%%%%%%%%%%%%%%

%%%%%%%%%%%%%%%%%%%%%%%%%%%%%%%%%%%%%%%%%%%%%%%%%%%%%%%%%%%%%%%%%%%%%%%%%%%%%%%
%%%%%%%%%%%%%%%%%%%%%%%%%%%%%%%%%%%%%%%%%%%%%%%%%%%%%%%%%%%%%%%%%%%%%%%%%%%%%%
\end{document}